\documentclass[pdflatex,sn-mathphys-num]{sn-jnl}% Math and Physical Sciences Numbered Reference Style

\usepackage{graphicx}%
\usepackage{multirow}%
\usepackage{amsmath,amssymb,amsfonts}%
\usepackage{amsthm}%
\usepackage{mathrsfs}%
\usepackage[title]{appendix}%
\usepackage{xcolor}%
\usepackage{textcomp}%
\usepackage{manyfoot}%
\usepackage{booktabs}%
\usepackage{algorithm}%
\usepackage{algorithmicx}%
\usepackage{algpseudocode}%
\usepackage{listings}%
\usepackage{tikz}%
\usepackage{array}
\usepackage{latexsym}
\def\ds{\displaystyle}

\mathchardef\isinpunto="0010

\newcommand{\abs}[1]{\left\vert#1\right\vert}

\def\R{I \!\! R}

\theoremstyle{thmstyleone}%
\newtheorem{theorem}{Theorem}%  meant for continuous numbers
\theoremstyle{thmstyletwo}%
\newtheorem{example}{Example}%

\theoremstyle{thmstylethree}%
\newtheorem{definition}{Definition}%

\newtheorem{nota}[theorem]{\bf Note}
\begin{document}

\title[P-curve: An improved solution for the single-facility location problem in two regions with $\ell_{p}$- and $\ell_{q}$-norms]
{P-curve: An improved solution for the single-facility location problem in two regions with $\ell_{p}$- and $\ell_{q}$-norms}


\author[1]{\fnm{Luis} \sur{Franco}}\email{lfranco@us.es}
\equalcont{These authors contributed equally to this work.}
\author[1]{\fnm{Francisco} \sur{Velasco}}\email{velasco@us.es}
\equalcont{These authors contributed equally to this work.}
\author*[1]{\fnm{F.J.} \sur{Ortega-Irizo}}\email{fjortega@us.es}
\author[1]{\fnm{Luis} \sur{Gonzalez-Abril}}\email{luisgon@us.es}
\equalcont{These authors contributed equally to this work.}

\affil*[1]{\orgdiv{Departamento de Econom\'ia Aplicada I}, \orgname{Universidad de Sevilla}, \orgaddress{\street{Avda. Ram\'on y Cajal,1}, \city{Sevilla}, \postcode{41018}, \country{Spain}}}


\abstract{The distance between two points is provided by the length of their test connecting path. However, there is no unique way to attain this path in the space $\R^2$, which is split by a straight line $L$ into two regions, $\Omega_p$ and $\Omega_q$, with $\ell_{p}$- and $\ell_q$-norms, respectively, where $1<p\leq q $ and $L\subset\Omega_q$.

In this paper, given a point $P\in \Omega_p$, the locus of points $Q\in\Omega_p$, called $p$-curve, is found, such that there are two ways to attain the distance between $P$ and $Q$. For each $P\in \Omega_p$, there are two $p$-curves which are two branches of two different $p$-parabola. These two branches split $\Omega_p$ into three subregions using a three-link shortest path in two of these regions. The implicit equation of each $p$-parabola is given.

Furthermore, by using the $p$-curve, an improved algorithm, called P-MFP, based on the MFP algorithm, is developed in order to solve the single-facility location problem in two regions with $\ell_{p}$ and $\ell_{q}$-norms. The results obtained with the new algorithm are better than those obtained with the MFP algorithm.ike subheadings, citations, or equations are permitted.}

%%=============



\maketitle
%%%%%%%%%%%%%%%%%%%%%%%%%%%%%%%%%%%%
\section{Introduction}%
%%%%%%%%%%%%%%%%%%%%%5%%%%%%%%%%%%%%
Let $L\equiv \, y= m \, x$ be a straight line which splits the space $\R^2$ into two regions, $\Omega_p$ and $\Omega_q$, with $\ell_{p}$- and
$\ell_q$-norms, respectively, where $1\leq p < q \leq \infty$ and
$L\subset\Omega_q$. Let us denote this construction as $\Omega_{pq}^m$.

Let $\cal{P}=$ $\{P_i\}_{i\in S_p}\subset \Omega_p$ and  $\cal{Q}=$$\{Q_j\}_{j\in S_q}\subset \Omega_q$ be two sets of points with $S_p$ and $S_q$ as two index sets, and  $w_i, w_j \in \R^{+}$ for any $i \in S_p,\, j\in S_q$. Hence, the following location problem is considered:
\begin{equation}\label{ec1-2}
\min_{X \in \R^2} \,D(X)=\min_{X \in \R^2} \, \left \{ \sum_{i\in S_p} w_i \, d(P_i,X)+
\sum_{j\in S_q} w_j \, d(Q_j,X)\right \}
\end{equation}

It is worth noting that real examples of the problem (\ref{ec1-2}) can be seen in \cite{Brimberg1} and \cite{Parlar}, where there are points in a city with streets drawn as horizontal and vertical lines on the map, and the rest of the points lie outside the city, and it is therefore possible to travel in a straight line. In \cite{Love2}, \cite{Blas} and \cite{Love}, two regions are also chosen, this time in London (The City and North London) and, by using the AD and SD criteria, a different $\ell_p$-norm in each region is obtained. Location problems, where a $\gamma$ norm with $\gamma= \ell_p$ or some block (or polyhedral) norm and $\phi$ is another block norm, is studied in \cite{Zafer-1}.

The problem (\ref{ec1-2}) was solved for $\ell_1$-norm, $\ell_2$-norm, and $L(x,y)\equiv  x= \alpha$ in \cite{Brimberg1} and  for $L\equiv y=m \, x$ in \cite{Parlar} and \cite{Zafer}. Furthermore,  \cite{Brimberg2} show a way to obtain the distance between certain points $X \in \Omega_p$ and $P \in \Omega_q$, based on three links. 
However, a more general study is carried out in \cite{Franco1}, where the shortest  paths with different $\ell_p$ norms are studied. These authors present a geodesic path connecting points $ P \in \Omega_p$ and $ Q \in \Omega_q$ and provide a point in $L$, called the ??Gate point??, by varying $p$- and $q$-distances, and the slope $m$ of the straight line $L$.

There are also papers, such as \cite{Blanco}, that use $\ell_p$ or polihedral (block) norms in $\R^n$ and consider only refraction: they do not use three-link shortest paths or a simple shortest path in a $d$-dimensional real space subdivided into several polyhedrals endowed with different $\ell_p$ norms  (\cite{Labbe,Gokhan}). Furthermore, \cite{Labbe} fail to take into account the properties developed in \cite{Franco2}, with norms $ \ell_1 $ and $ \ell_q $, where it is shown that, depending on the localisation of the points, the distance is obtained from a path of one, two, or three links. That is to say, for three links there are points $ P \in \Omega_1 $ and $ Q \in \Omega_1 $ for which it is possible to obtain two points $ P_1, Q_1 \in L$ such that $ d_1 (P, P_1) + d_q (P_1, Q_1) + d_1 (Q_1, Q) \leq d_1 (P, Q) $. 

The Location problem (\ref{ec1-2}) is very well studied in \cite{Plastria-I,Plastria-II} with $card(S_p)= card(S_q)=1$, where the single, double (Sections 4.2-4.3 \cite{Plastria-I}) and triple-link path (Section 4.4 \cite{Plastria-I}) are studied in detail,  and where $\ell_p$ and $\ell_q$ norms are generalised by $\phi$ and $\gamma$ gauges in the space $\R^d$. The paper (\cite{Plastria-II}) is dedicated to the particular case in which $ d = 2 $ and a semi-$\phi$-parabola is found as the locus formed by the set of points where $d_{p}(P,Q_1)+d_{q}(P_1,Q_1)+d_{p}(Q_1,Q) = d_p(P,Q)$, with $P_1$ and $Q_1$ as described above. Figure 5 in \cite{Plastria-II} shows that,  given a $P \in L_{\phi}$ and a $w$ parallel to the $L$ line, with $\phi(w) > \gamma(w)$, a point $X$ is constructed that belongs to a curve that represents a semi $\phi$-parabola,  In Figure 6 in \cite{Plastria-II}, the semi -$ \phi $-parabola is given, where several points X of the $\phi-$ parabola are obtained. Consequently, if a point $Q \in L_{\phi}$ is above the semi $\phi$-parabola and below the $L$ line, then a three link path is necessary to obtain the distance defined in \cite{Plastria-I} between $P$ and $Q$, where $L_{\phi}=\Omega_{\phi}$. 
Therefore, in order to solve the problem (\ref{ec1-2}), it is important to take into account the distance $d(\cdot,\cdot)$ studied in \cite{Brimberg1} and \cite{Plastria-I}. 

The main aim of this paper is twofold: firstly,  to find an analytical expression of the locus of the points in $\Omega_{pq}^m$ such that there are two different ways (single or triple-link path) of reaching a given point $P$ from $Q$, where $P$ and $Q$ are in $\Omega_p$; and secondly, to solve the problem (\ref{ec1-2}) by using the result above. To this end,  a procedure, called P-MFP, is developed which improves the MFP procedure (\cite{Brimberg1}). Hence these papers can be improved by taking into account the results obtained in this paper. 

The remainder of this paper is arranged as follows:
the distance between two points in $\Omega_{pq}^m$\ is studied in Section~\ref{sec:distance}.
Section~\ref{sec:pcurve} presents a theoretical study of locus of the points in $\Omega_{pq}^m$\ such that there are two different ways of reaching a given point $P$, and implicit equations of two $p$-parabolas are found that generalise the paper of \cite{Franco2}.
Based on the previous results, the new procedure; P-MFP, is developed in order to solve the problem  (\ref{ec1-2}) in Section~\ref{section_pMFP} and an experimentation is carried out in Section~\ref{section6}.
Lastly, conclusions are drawn in Section 6.

%%%%%%%%%%%%%%%%%%%%%%%%%%%%%%%%%555
\section{Distance between two points in $\Omega_{pq}^m$}\label{sec:distance}
%%%%%%%%%%%%%%%%%%%%%%%%%%%%%%%%%

\begin{definition}
Given two points $P,\, Q\in \R^2$, the distance between them is defined as the length of the
shortest path (geodesic path) connecting $P$ and $Q$ (\cite{Mitchell}), denoted by $d(P,Q)$.
\end{definition}

\begin{figure}[ht]\centering
\begin{tikzpicture}[xscale=0.5,yscale=0.5]
\draw[-,thick] (0,0) --(7,7)node[below]{\footnotesize $\phantom{xxxxxxx}L\equiv y= x$};
%\draw[ultra thin, help lines] (0,0) grid (10,5);
\node at (4,1) {\large $\Omega_{1}$};
\node at (1,4) {\large $\Omega_2$};
\draw[red,thick] (1,0)--(1,1)--(5,5)--(6,5);
\draw[blue,thick] (1,0)--(6,0)--(6,5);
\filldraw [black] (1,0) circle (2.5pt) node[below] {\footnotesize $P_1(1,0)$};
\filldraw [black] (6,5) circle (2.5pt) node[right] {\footnotesize $P_2(6,5)$};
\filldraw [black] (1,1) circle (2.5pt) node[left] {\footnotesize $Q_1(1,1)$};
\filldraw [black] (5,5) circle (2.5pt) node[left] {\footnotesize $Q_2(5,5)$};
\end{tikzpicture}
\begin{tikzpicture}[xscale=0.5,yscale=0.5]
\draw[-,thick] (0,0) --(7,7)node[below]{\footnotesize $\phantom{xxxxxxx}L\equiv y= x$};
%\draw[ultra thin, help lines] (0,0) grid (10,5);
\node at (4,1) {\large $\Omega_2$};
\node at (1,4) {\large $\Omega_3$};
\draw[red,thick] (1,0)--(1,1)--(5,5)--(6,5);
\draw[blue,thick] (1,0)--(6,5);
\filldraw [black] (1,0) circle (2.5pt) node[below] {\footnotesize $P_1(1,0)$};
\filldraw [black] (6,5) circle (2.5pt) node[right] {\footnotesize $P_2(6,5)$};
\filldraw [black] (1,1) circle (2.5pt) node[left] {\footnotesize $Q_1(1,1)$};
\filldraw [black] (5,5) circle (2.5pt) node[left] {\footnotesize $Q_2(5,5)$};
\end{tikzpicture}
\caption {Alternative trajectories (red and blue) between two points $P_1$ and $P_2$ in $\Omega_{p\,q}^m$.}
 \label{fig:transyectoria}
\end{figure}
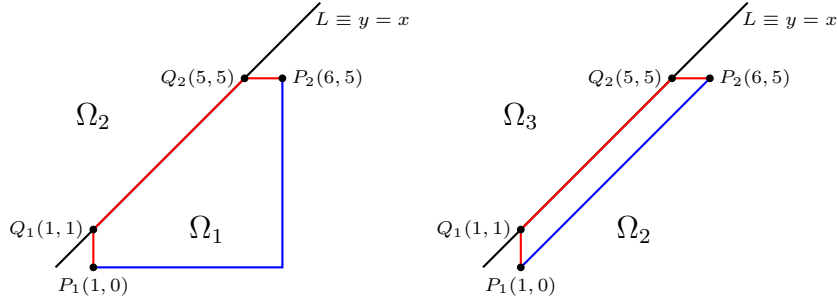

Let $P,\, Q\in \Omega_{pq}^m$ be two points, and $d_p(\cdot,\cdot)$ and $d_q(\cdot,\cdot)$ be the induced distances from $\ell_p$- and $\ell_q$-norms, respectively. It is trivial that if $p=q$ then $d(P,Q)=d_p(P,Q)=d_q(P,Q)$: that is, the straight line $L$ is not relevant. Nevertheless,  if $p\neq q$, then the shortest path in $\Omega_{pq}^m$ is not always given by the segment joining $P,Q \in \Omega_p$ (\cite{Mitchell}, \cite{Franco2}, \cite{Plastria-I}, \cite{Plastria-II}) and, in these cases, the straight line $L$ is relevant. Let us study a pair of  counterexamples.
\begin{itemize}
\item Let $\Omega_{1\,2}^1$ be a region, where the points $P_1(1,0)$ and $P_2(6,5)\in\Omega_1$;  $Q_1(1,1)$ and $Q_2(5,5)\in L\subset\Omega_2$.

Hence,  by definition (11) given in \cite{Plastria-I}: $d(P_1,P_2) \leq d(P_1,Q_1)+d(Q_1,Q_2)+ d(Q_2,P_2)=d_1(P_1,Q_1)+d_2(Q_1,Q_2)+ d_1(Q_2, P_2)=   1+4\,\sqrt{2}+1 < 5+5= d_1(P_1,P_2)$ (see \figurename{ \ref{fig:transyectoria}} left-hand-side).

\item Let $\Omega_{2\,3}^1$ be a region, where the points $P_1(1,0)$ and $P_2(6,5)\in\Omega_2$;  $Q_1(1,1)$ and $Q_2(5,5)\in L\subset\Omega_3$.

Hence,  by definition (11) given in \cite{Plastria-I}: $d(P_1,P_2) \leq d(P_1,Q_1)+d(Q_1,Q_2)+ d(Q_2,P_2)=d_2(P_1,Q_1)+d_3(Q_1,Q_2)+ d_2(Q_2, P_2)=
  1+4\,\sqrt[3]{2}+1 < 5\sqrt{2}= d_2(P_1,P_2)$ (see \figurename{ \ref{fig:transyectoria}} right-hand-side).
\end{itemize}
%%%%%%%%%%
% %%%%%%%%%%%%%%%%%%%%%%%%%%%%%%%%%%%%%%%%%%%%%%%%%%%%%%%%%%

Therefore,  the location of the points in $\Omega_{pq}^m$ are of major importance for the distance between them to be obtained. In order to further elaborate on the subject, let us study the three possibilities depending on  the location of the points:
\begin{enumerate}
 \item If $P,\,Q\in \Omega_q$ then it is straightforward to prove that $d(P,Q)=d_q(P,Q)$ for the triangle inequality in $\Omega_q$ and  that  $d_q(\cdot,\cdot)\leq d_p(\cdot,\cdot)$ because $p<q$.
\item If $P\in \Omega_p$ and $Q\in \Omega_q$, then the length of the shortest path connecting $P$ and $Q$ is given as follows:
\begin{equation}\label{ec1-1}
d(P,Q)=\min_{X \in L} \, \left \{d_p(P,X)+d_q(X,Q)\right \}.
\end{equation}
\cite{Brimberg1}, show that the shortest path between two points $P$ and $Q$ in regions with $\ell_1$ and
$\ell_2$ norms, respectively, is the polygonal $\overline{PM} \cup \overline{MQ}$, where $M$ is the orthogonal projection of $P$ on the line separating the regions, when considering only vertical and horizontal lines. \cite{Zafer} generalise the result of any straight line and norms $\ell_1$ and $\ell_2$.\\
 However, the  problem  (\ref{ec1-1}) has been studied and explicitly solved in \cite{Franco1}, and also studied in \cite{Plastria-I} for gauges $\phi$, $\gamma$. Hence, it is always possible to find a point $G\equiv G(P,Q)\in L$, called the Gate point, of $P$ and $Q$ (\cite{Brimberg2}), such that $$d(P,Q)=d_p(P,G)+d_q(G,Q),\qquad \mbox{if } P\in\Omega_p\mbox{ and } Q\in \Omega_q. $$
 An example can be observed in \figurename{ \ref{fig:Gate}}.
 %%%%%%%%%%
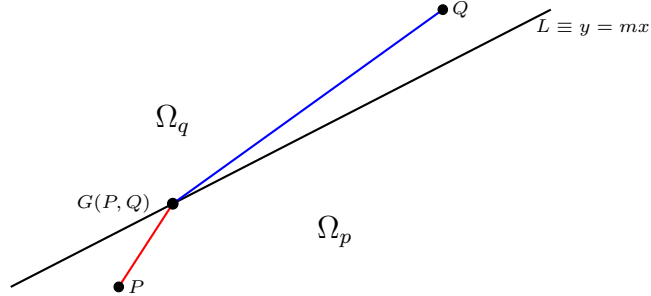
\begin{figure}[ht]\centering
\begin{tikzpicture}[xscale=1.00/1.4,yscale=1.10/1.5]
\draw[-,thick] (0,0) --(10,5)node[below]{\footnotesize $\phantom{xxxxxxx}L\equiv y=m x$};
%\draw[ultra thin, help lines] (0,0) grid (10,5);
\node at (6,1) {\large $\Omega_{p}$};
\node at (3,3) {\large $\Omega_q$};
\draw[red,thick] (2,0)-- (3,1.5);
\draw[blue,thick] (8,5)-- (3,1.5);
\filldraw [black] (2,0) circle (2.5pt) node[right] {\footnotesize $P$};
\filldraw [black] (8,5) circle (2.5pt) node[right] {\footnotesize $Q$};
\filldraw [black,thick] (3,1.5) circle (2.5pt) node[left] {\footnotesize $G(P,Q)\phantom{x}$};
\end{tikzpicture}
\caption {Geodesic path between $P\in\Omega_p$ and $Q\in\Omega_q$. The $G(P,,Q)$ point is the Gate point of $P$ and $Q$.}
 \label{fig:Gate}
\end{figure}
% %%%%%%%%%%%%%%%%%%%%%%%%%%%%%%%%%%%%%%%%%%%%%%%%%%%%%%%%%%
It is worth noting that the problem (\ref{ec1-1}) can be seen as the simplest location problem in two different regions. 

To better understand how the ``Gate''  is obtained, it is necessary to expalin how the distance is accertained.

First let us remember that the distance  $p$ with $1<p$ from a point $P= (a,b)$ to a straight line $L(x,y)\equiv \omega_1 \, x+ \omega_2 \, y + \omega_0=0$, denoted by $d_p(P,L)$,  is:
$$d_p(P,L)=\frac{\abs{ \omega_1 \, a +\omega_2\, b + \omega_0}}{|\omega_1|^{\bar{p}}+|\omega_2|^{\bar{p}}} $$
where $\bar{p}$ is the conjugate of $p$, that is, $\bar{p}=\frac{p-1}{p}$.

The straight line with $\omega_2=0$ is studied in \cite{Brimberg1}, and it can henceforth be considered that $\omega_2\neq 0$ and that it is always possible to select the straight line as $L\equiv y= m \, x$ (\cite{Franco1}). 
By taking  the line $L\equiv y= m \, x$ into consideration, the following definition is given in \cite{Franco2}:
\begin{definition} (Projection)
Let $X= (x_0,y_0)$ be a point and $L \equiv y= m \, x$ be a straight line in $\R^2$. The projection point $X^*$ of $X$ on the straight line $L$ with respect to direction $m^*$, with $m\neq m^*$  (see \figurename{ \ref{fig:projection}}), is defined as follows:
$$X^* = Proj(X,L, m^*) = (x^*, y^*) = \left (\frac{m^*\, x_0- y_0}{m^* - m},\frac{m\,(m^*\, x_0- y_0)}{m^* - m} \right )$$
\end{definition}
%%%%%%%%%%
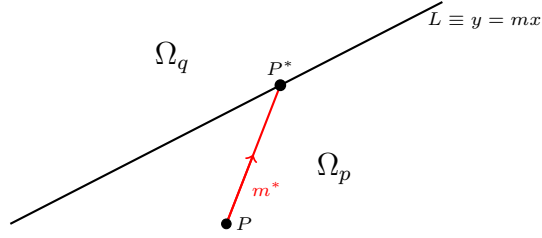
\begin{figure}[ht]\centering
\begin{tikzpicture}[xscale=1.0/1.4,yscale=1.10/1.5]
\draw[-,thick] (0,0) --(8,4)node[below]{\footnotesize $\phantom{xxxxxxx}L\equiv y=m x$};
%\draw[ultra thin, help lines] (0,0) grid (10,5);
\node at (6,1) {\large $\Omega_{p}$};
\node at (3,3) {\large $\Omega_q$};
\draw[red,thick,->] (4,0)-- (4.5,1.25);
\draw[red,thick] (4,0)-- (5,2.5);
\node[red] at (4.75,0.625) {\footnotesize$m^*$};
\filldraw [black] (4,0) circle (2.5pt) node[right] {\footnotesize $P$};
\filldraw [black,thick] (5,2.5) circle (2.5pt) node[above] {\footnotesize $P^*$};
\end{tikzpicture}
\caption {The projection point $P^*$ of $P$ on the straight line $L$ with respect to direction $m^*$.}
 \label{fig:projection}
\end{figure}
% %%%%%%%%%%%%%%%%%%%%%%%%%%%%%%%%%%%%%%%%%%%%%%%%%%%%%%%%%%%%%%%%%%
 
 Let $P=(a,b)\in \Omega_p$ and $Q=(c,d) \in L \subset \Omega_q$ be points, $1 < p< q < \infty$, in the plane separated by the straight line $L\equiv y= m \, x$ with $m \geq 0$. Suppose further that $G= G(P,Q)$ is the Gate point, and therefore either $G(P,Q)=Q$ or $G(P,Q)=Proj(P,L, m_p)$, where $m_p$ is one of the solutions of the equation (\ref{precta}) (see \figurename{ \ref{fig:two-link}}):
 \begin{equation}\label{precta}
 \frac{{1 + |m| \cdot \left| m_p \right|^{p - 1} }}{{\left( {\sqrt[p]{{1 + |m_p|^p }}} \right)^{p - 1} }} = \sqrt[q]{{1 + m^q }}
\end{equation}
%%%%%%%%%%
\begin{figure}[ht]\centering
\begin{tikzpicture}[xscale=1.0/1.4,yscale=1.10/1.5]
\draw[-,thick] (0,0) --(10,5)node[below]{\footnotesize $\phantom{xxxxxxx}L\equiv y=m x$};
%\draw[ultra thin, help lines] (0,0) grid (10,5);
\node at (6,1) {\large $\Omega_{p}$};
\node at (3,3) {\large $\Omega_q$};
\draw[red,thick] (2,0)-- (3,1.5);
\draw[blue,thick] (8,4)-- (2,0);
\filldraw [black] (2,0) circle (1.5pt) node[right] {\footnotesize $P$};
\filldraw [black] (8,4) circle (1.5pt) node[above] {\footnotesize $Q$};
\filldraw [black,thick] (3,1.5) circle (1.5pt) node[left] {\footnotesize $Q=G(P,Q)\phantom{x}$};
\end{tikzpicture}
\caption {The  red and blue straight lines represent the two possible paths between $P$ and $Q\in L $.}
 \label{fig:two-link}
\end{figure}
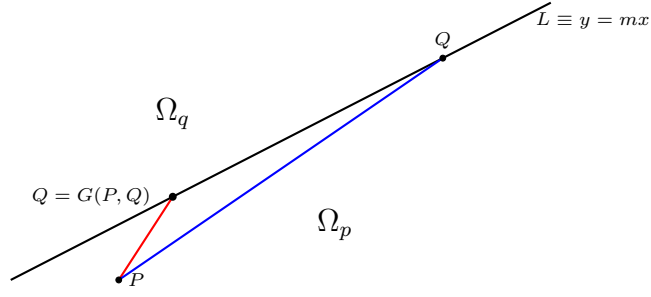
% %%%%%%%%%%%%%%%%%%%%%%%%%%%%%%%%%%%%%%%%%%%%%%%%%%%%%%%%%%
The proof can be found in \cite{Franco0}.
Thus, for example, if $p=2$, then an algebraic expression of the solution is as follows: $$m_2= \frac{m-k \, \sqrt{m^2+1-k^2}}{k^2-m^2}<\frac{m+k \, \sqrt{m^2+1-k^2}}{k^2-m^2}=m_1,$$ where $k= \sqrt[q]{1+m^q}$.

 \item If $P,\,Q\in  \Omega_p$ then, by taking into account that $d_q(\cdot,\cdot)\leq d_p(\cdot,\cdot)$,  there are two possibilities of obtaining $d(P,Q)$: $d_p(P,Q)$ and $$\min_{X_1, X_2\in L} \left\{ d_p(P,X_1)+ d_q(X_1,X_2)+d_p(X_2,Q) \right\}$$ (see \figurename{ \ref{fig:rand1}}). This property is shown in \cite{Franco2} for $p=1$, $p\leq q$ and  generalised in \cite{Plastria-I} (Corollary 4) for gauges $\phi$, $\gamma$.
\end{enumerate}
 % %%%%%%%%%%%%%%%%%%%%%%%%%%%%%%%%%%%%%%%%%%%%%%%%%%%%%%%%%%
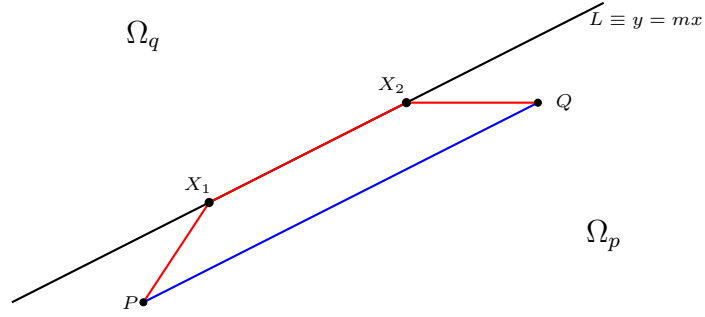
\begin{figure}[ht]\centering
\begin{tikzpicture}[xscale=1.0/1.15,yscale=1.10/1.25]
\draw[-,thick] (0,0) --(9,4.5)node[below]{\footnotesize $\phantom{xxxxxxx}L\equiv y=m x$};
%\draw[ultra thin, help lines] (0,0) grid (10,5);
\node at (9,1) {\large $\Omega_p$};
\node at (2,4) {\large $\Omega_q$};
\draw[red,thick] (3,1.5)-- (2,0);
\draw[blue,thick] (2,0)--(8,3);
\draw[red,thick] (3,1.5)-- (6,3)--(8,3);
\node at (8.4,3){\footnotesize $Q$};
\node at (1.8,0){\footnotesize $P$};
\filldraw [black,thick] (3,1.5) circle (1.5pt) node[above] {\footnotesize $X_1$\phantom{xx}};
\filldraw [black,thick] (6,3) circle (1.5pt) node[above] {\footnotesize $X_2$ \phantom{xx}};
\filldraw [black] (8,3) circle (1.5pt);
\filldraw [black] (2,0) circle (1.5pt);
\end{tikzpicture}
\caption {The  red and blue straight lines represent the two possible paths between $P$ and $Q$.}
 \label{fig:rand1}
\end{figure}
% %%%%%%%%%%%%%%%%%%%%%%%%%%%%%%%%%%%%%%%%%%%%%%%%%%%%%%%%%%

%%%%%%%%%%%%%%%%%%%%%%%%%%
\section{The $p$-curve}\label{sec:pcurve}

The localisation problem (\ref{ec1-2}) was solved in \cite{Franco2}  with $p=1$ and $q>1$. When studying the localisation problem in $\Omega_1$, they find three regions defined by a horizontal line and a vertical line, which they call $\Omega_{1R}, \, \Omega_{1C}, \mbox{ and}\, \Omega_{1L}$. These lines are the locus of points equidistant from a point $P=(a,b) \in \Omega_{1C}$ and a point $Q=(c,d) \in \Omega_{1R}$ or $\Omega_{1L}$, respectively (\figurename{} 4 in \cite{Franco2}).
In this section, the loci for $p>1$ is studied.

Let $\Omega _{pq}^m $ be with $1 < p< q$, and a fixed point $A=(a,b)\in  \Omega_p$. In order to calculate $d(A,P)$ for $P \, \in  \Omega_p$, it is necessary to obtain $d_p(A,P)$ (see Section~\ref{sec:distance})  and 
\begin{equation}\label{eq_dospuntos}
\min \left\{ d_p(A,X_1)+d_q(X_1,X_2)+d_p(X_2,P);\,\forall X_1,X_2\in L\right  \}.
\end{equation}

Hence, the aim is identify  the locus of points $P=(x,y) \in \Omega_p$, denoted by $\Delta$, such that
\begin{equation}\label{Delta}
\Delta\equiv \left \{P\in \Omega_p \mbox{ such that } d_p (A,A') + d_q (A',P') + d_p (P',P) - d_p (A,P)=0
\right \}
\end{equation}
where $A'$ and $P'$ are the points $X_1$ and $X_2$, which are solution of the problem (\ref{eq_dospuntos}).
%$d_p (A,A') + d_q (A',P') + d_p (P',P)=\min\{d_p(A,X_1)+d_q(X_1,X_2)+d_p(X_2,P);\,\forall X_1,X_2\in L\}$.

To this end, it is necessary to obtain the points $A',P'\in L$ which are identified as the Gate points (\cite{Franco1})   from the equation introduced in Note 2.12 of \cite{Franco1}:
\begin{equation}\label{fr5}
\frac{{1 + m \cdot \left| {m'} \right|^{p - 1} }}{{\left( {\sqrt[p]{{1 + \left| {m'} \right|^p }}} \right)^{p - 1} }} 
= \sqrt[q]{1 + m^q }.
\end{equation}
In Proposition 2.22 of \cite{Franco0}, it is proved that this equation has two solutions: the  $m_1$, $m_2$ functions of $m$, such that $0< m_2 < m< m_1$.

In order to find $\Delta$, there are two cases: $a<x$ and  $a>x$. These provide two  branches on $\Delta$, called the Right and Left  branches on $\Delta$, respectively.

%%%%%%%%%%%%%%%%%%%%%%%%%%%%%%%%%%%%%%%555
\subsection{Right branch on $\Delta$ ($\Delta_r$) }
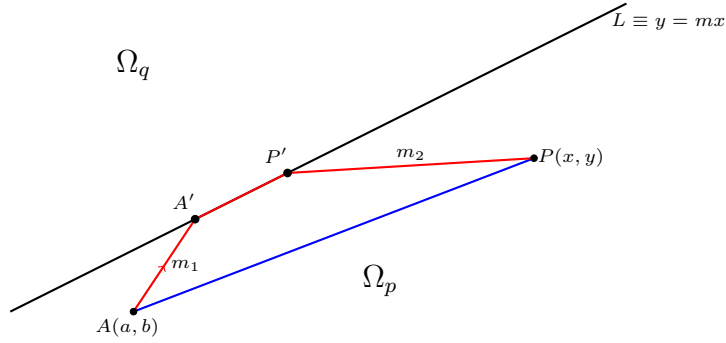
\begin{figure}[ht]\centering
\begin{tikzpicture}[xscale=1.10/1.35,yscale=1.10/1.35]
%\draw[fill=gray!15,ultra thin,dotted] (3,1.5)--(10,1.86)--(10,5)--(3,1.5) ;
\draw[-,thick] (0,0) --(10,5)node[below]{\footnotesize $\phantom{xxxxxxx}L\equiv y=m x$};
%\draw[ultra thin, help lines] (0,0) grid (10,5);
\node at (6,0.5) {\large $\Omega_p$};
\node at (2,4) {\large $\Omega_q$};
\node at (1.9,-0.25){\footnotesize $A(a,b)$};
\draw[red,thick] (3,1.5)-- (2,0);
\draw[blue,thick] (2,0)--(8.5,2.49);
\draw[red,thick] (3,1.5)-- (4.5,2.25)--(8.5,2.49);
\filldraw [black] (8.5,2.49) circle (1.5pt);
\node at (9.1,2.5){\footnotesize $P(x,y)$};
\node at (2.85,0.75){\footnotesize $m_1$};
\node at (6.5,2.55){\footnotesize $m_2$};  
\draw[->,red] (2,0)--(2.5,.75);
\filldraw [black] (2,0) circle (1.5pt);
\filldraw [black,thick] (4.5,2.25) circle (1.5pt) node[above] {\footnotesize $P'$\phantom{xx}};
\filldraw [black,thick] (3,1.5) circle (1.5pt) node[above] {\footnotesize $A'$\phantom{xx}};
%\draw[thick,dashed](3,1.5)--(10,1.86);
\end{tikzpicture}
\caption {Situation to obtain the Right branch on $\Delta$.}
 \label{fig:ladoderecho}
\end{figure}

By taking into account that $a<x$, then the Gate points $A'= Proj(A,L, m_1) $ and $P'=Proj(P,L, m_2)$, that is,
$$A'=\left( {\frac{{b - m_1  \cdot a}}{{m - m_1 }}  ,m\frac{{b - m_1  \cdot a}}{{m - m_1 }}} \right)\quad\mbox{  and  }\quad
P'=\left(  {\frac{{y - m_2  \cdot x}}{{m - m_2 }},m\frac{{y - m_2  \cdot x}}{{m - m_2 }}} \right).$$
Furthermore, in order to have a three-link path, the abscissa of point $A'$ must be less than the abscissa of point $P'$, that is, 
\begin{equation}\label{conditionpara}
\frac{b - m_1  \cdot a}{m - m_1 } < \frac{y - m_2  \cdot x}{m - m_2 }
\end{equation}
 (this condition is depicted in grey in \figurename{ \ref{fig:ladoderecho}}).
%%%%%%%%%%%%%%%%%%%%%
% %%%%%%%%%%%%%%%%%%%%%%%%%%%%%%%%%%%%%%%%%%%%%%%%%%%%%%%%%%%%%%

Hence, the analytical expression obtained  from (\ref{Delta}) of the Right branch on  $\Delta$ ($a< x$),  denoted by $\Delta_r$, is as follows:
\begin{equation}\label{ec1}
\Delta_r  \equiv \left\{ P\in \Omega_p \mbox{ such that }A_r x + B_r y + C_r - \sqrt[p]{{(x - a)^p  + (y - b)^p }}=0\right\}
\end{equation}
 where %it has been taken into account the equation (\ref{fr2}) obtained from (\ref{fr5}) \cite{Franco1}:\\
 \begin{equation}\label{ec1-3}
 A_r = \ds\frac{1}{{\left( {\sqrt[p]{{1 + m_2^p }}} \right)^{p - 1} }},
 B_r =\ds \left( {\frac{{m_2 }}{{\sqrt[p]{{1 + m_2^p }}}}} \right)^{p - 1}, \mbox{ and }
C_r = \ds - \frac{{a + m_1^{p - 1} b}}{{\left( {\sqrt[p]{{1 + m_1^p }}} \right)^{p - 1} }}
\end{equation}
and
 \begin{equation}
\frac{{1 + m \cdot m_2^{p - 1} }}{{\left( {\sqrt[p]{{1 + m_2^p }}} \right)^{p - 1} }} = \frac{{1 + m \cdot m_1^{p - 1} }}{{\left( {\sqrt[p]{{1 + m_1^p }}} \right)^{p - 1} }} = \sqrt[q]{{1 + m^q }} \label{fr2}
\end{equation}

Hence, let us show some properties:
\begin{itemize}
\item $ 0 < A_r , B_r < 1 $

\item  $A_r x + B_r y + C_r= d_p (A,A') + d_q (A',P') + d_p (P',P)$ by (\ref{ec1}), that is, $A_r x + B_r y + C_r$ provides the distance of the three-link path that joins point $A$ with point $P$ passing through points $A'$ and $P'$ on the line $L$.

\item By considering $\bar{p}$, the conjugate of p ($\bar{p}= \frac{p}{p-1}$), it follows that:
$$A_r^{\bar{p}}+B_r^{\bar{p}}=\displaystyle{\frac{1}{(1+m_2^p)^{\frac{p-1}{p}\frac{p}{p-1}}}+\frac{m_2^{(p-1)\frac{p}{p-1}}}{(1+m_2^p)^{\frac{p-1}{p}\frac{p}{p-1}}}=\frac{1+m_2^p}{1+m_2^p}=1}$$

\item It is worth noting that by setting $1<p<q$ and $A \in \, \Omega _p $, the $\Delta_r $ provides a curve which will be called the   $p$-curve. Thus, $\Delta_r $ splits  the half-plane $\Omega _p $ in such a way that:
\begin{itemize}
\item if $P$ is above the curve,  then 
$\Delta_r<0 \equiv A_r x + B_r y + C_r < \sqrt[p]{{(x - a)^p  + (y - b)^p }},
$ and the shortest path between point $A$ and  point $P(x,y)$ is a polygonal with three segments with distance $A_r x + B_r y + C_r$. Let us denote this region as $R_{pr}$ (see \figurename{ \ref{fig:slopes}}).    
\item if $\Delta_r> 0$, then the shortest path between point $A$ and  point $P(x,y)$ is the segment joining these points with distance $d_p(A,P)$. 
\item otherwise, $\Delta_r= 0$, and the two paths between points $A$ and $P(x,y)$ have the same distance.
\end{itemize}
\end{itemize}

The definition of parabola is: ``Geometric locus of points $P(x,y)$ equidistant  from a straight line (directrix) and a fixed point (focus)''.  Hence, by taking into account the $\ell_p$-distance, it considers the parabola, called $p$-parabola and denoted ${\cal P}_r$, of the directrix $L_r\equiv A_r x + B_r y + C_r=0$ and the focus point $A(a,b)$. 

Let us see several features of the $p$-parabola ${\cal P}_r$.
\begin{itemize}
 \item The vertex of the $p$-parabola ${\cal P}_r$ is $\ds
V_r=\left(\frac{a+x_0}{2}\,,\, \frac{b+y_0}{2}\right )$
where
$$
x_0= \ds \frac{a\, B_r^2- A_r (b \, B_r + C_r)}{A_r^2+ B_r^2}, \mbox{ and }
y_0= \frac{b\, A_r^2- B_r (a \, A_r + C_r)}{A_r^2+ B_r^2}
$$
obtained after solving the system:
\begin{equation}\label{ec1-4}
\left\{\begin{array}{rcl}
A_r x+ B_r y +C_r& = & 0\\
-B_r x + A_r y + B_r a -A_r b&=& 0
\end{array}\right .
\end{equation}
where the second equation %of the system (\ref{ec1-4}) 
is the axis of the $p$-parabola ${\cal P}_r$.
 %(See Fig.\ref{parabolas}):\\

\item Point $A'= Proj(A,L, m_1)=\left( {\frac{{b - m_1  \cdot a}}{{m - m_1 }}  ,m\frac{{b - m_1  \cdot a}}{{m - m_1 }}} \right) $ belongs to  ${\cal P}_r$.\\
{\it Proof: }\\
By using equations (\ref{ec1-3}) and (\ref{fr2}), it follows that:\\
$$\begin{array}{lcr}
A_r  \left(\frac{b-m_1\, a}{m-m_1}\right)+ B_r \left(m\,\frac{b-m_1\, a}{m-m_1}\right)+C_r-d_p(A,A')&=&\cdots\\  
\left ( \frac{b-m_1\, a}{m-m_1}\right )\left (\frac{1+m\, m_2^{p-1}}{(\sqrt[p]{1+m_2^p}\,)^{p-1}} \right )
-\frac{a+ m_1 \, b}{\left(\sqrt[p]{1+m_1^p}\,\right)^{p-1}}-
\frac{b- m \, a}{m-m_1}\, \sqrt[p]{1+m_1^p}&=&\cdots\\
\mbox{By taking into account  the first equality in equation (\ref{fr2})}& &\\
\left (\frac{b-m_1\, a}{m-m_1}\right )\, \left ( \frac{1+m\, m_1^{p-1}}{\left(\sqrt[p]{1+m_1^p}\,\right)^{p-1}}\right )
-\frac{(a+m_1^{p-1}\,b)\,(m-m_1)-(b- m \, a)\,(1+m_1^p)}{(m-m_1)\left(\sqrt[p]{1+m_1^p}\,\right)^{p-1} }&=&\cdots\\
 \frac{ (b-m_1 a)(1+ m \, m_1^{p-1})-
 (a+ m_1^{p-1}  b)(m-m_1) - (b- m a)(1+m_1^p)}{(m-m_1) \left(\sqrt[p]{1+m_1^p}\right)^{p-1}} &=&\cdots\\
 \frac{b+b\, m\,m_1^{p-1}-m_1\,a+a\,m\,m_1^p-(a\, m-a \,m_1+m\, m_1^{p-1} b-m_1^p b)-(b+b\, m_1^p- m\, a-a\, m \,m_1^p)}{(m-m_1) \,\left(\sqrt[p]{1+m_1^p}\,\right)^{p-1}} \, &=&0
\end{array}$$
  
As a consequence of the previous result, if $P$ satisfies condition  (\ref{conditionpara}), then $A_r x + B_r y + C_r>0$ (since it is a distance) and, since
$$d_p(P, L_r)=\frac{\abs{A_r x + B_r y + C_r}}{\sqrt[\bar{p}]{A_r^{\bar{p}}+B_r^{\bar{p}}}}=
A_r x + B_r y + C_r$$
it follows that  $\Delta_r $ represents, from point $A'=\left( {\frac{{b - m_1  \cdot a}}{{m - m_1 }}  ,m\frac{{b - m_1  \cdot a}}{{m - m_1 }}} \right)$, the Right branch of the $p$-parabola ${\cal P}_r$ (see \figurename{ \ref{fig:slopes}}).

\item The tangent line  to  $p$-parabola ${\cal P}_r$ at point $A'$ has the slope of the straight line $L\equiv y= m \, x$.\\
{\it Proof: }\\
If $A'$ is denoted as $A'=(a_0,b_0)$, then 
$$d_p(A,A')=  \sqrt[p]{{(a_0 - a)^p  + (b_0 - b)^p }}=(a_0-a) \sqrt[p]{1+ m_1^{p}}$$
Hence: ${\cal P}_r\equiv A_r x+B_r\,y+ C_r-d_p(A,P)=0$ and 
the slope of the tangent line to ${\cal P}_r$ in $A'$ is:\\
$$\begin{array}{lcl}
-\ds\frac{\frac{\partial{\cal P}_r (a_0,b_0)}{\partial x}}{\frac{\partial {\cal P}_r (a_0,b_0)}{\partial y}}&=&
\frac{-A_r+\frac{(a_0-a)^{p-1}}{\left(\sqrt[p]{{(a_0 - a)^p  + (y_o - b)^p }}\,\right )^{p-1}}}{B_r-\frac{(b_0-b)^{p-1}}{\left (\sqrt[p]{{(a_0 - a)^p  - (b_0 - b)^p }}\,\right )^{p-1}}}\\[.45cm]  
&=&\frac{-A_r\,(a_0-a)^{p-1} \left(\sqrt[p]{{1  + m_1^p }}\,\right)^{p-1}+(a_0-a)^{p-1}}{B_r\,(a_0-a)^{p-1} \left(\sqrt[p]{{1+m_1^p }}\, \right )^{p-1}-m_1^{p-1}(a_0-a)^{p-1}}\\[.45cm]
  
&=&
\frac{-A_r \left(\sqrt[p]{{(a_0 - a)^p  + (b_0 - b)^p }}\,\right)^{p-1}+(x-a)^{p-1}}{B_r \left(\sqrt[p]{{(a_0 - a)^p  - (b_0 - b)^p }}\,\right )^{p-1}-(b_0-b)^{p-1}}\\[.45cm]

&=&\frac{-A_r\, \left (\sqrt[p]{{1  + m_1^p }}\,\right)^{p-1}+1}{B_r\, \left(\sqrt[p]{{1+m_1^p }}\, \right)^{p-1}-m_1^{p-1}}\\[.45cm]

&=&\frac{- \frac{\left(\sqrt[p]{{1  + m_1^p }}\,\right)^{p-1}}{\left(\sqrt[p]{{1  + m_2^p }}\,\right)^{p-1}}+1}
{\frac{m_2^{p-1}\,\left(\sqrt[p]{{1  + m_1^p }}\,\right)^{p-1}}{\left(\sqrt[p]{{1  + m_2^p }}\,\right)^{p-1}}-m_1^{p-1}
}\\[.65cm]
  
&=&\frac{\sqrt[p]{{1  + m_2^p }}\,)^{p-1}-(\sqrt[p]{{1  + m_1^p }}\,)^{p-1}}{m_2^{p-1}\,(\sqrt[p]{{1  + m_2^p }}\,)^{p-1}-m_1^{p-1}(\sqrt[p]{{1  + m_1^p }}\,)^{p-1}}\\[.45cm]

&=&\frac{\left(\sqrt[p]{{1  + m_2^p }}\,\right)^{p-1}-\left(\sqrt[p]{{1  + m_1^p }}\,\right)^{p-1}}{m_2^{p-1}\,\left(\sqrt[p]{{1  + m_2^p }}\,\right)^{p-1}-m_1^{p-1}\left(\sqrt[p]{{1  + m_1^p }}\,\right)^{p-1}}  
\end{array}$$
Taking into account the first equality of equation (\ref{fr2}), it follows that:\\
$$\left(1+m \, m_2^{p-1}\right)\left(\sqrt[p]{1  + m_1^p }\,\right)^{p-1}=(1+m \, m_1^{p-1})\left(\sqrt[p]{1  + m_2^p }\,\right)^{p-1}$$  
 and hence
$$ m \, m_2^{p-1} \left(\sqrt[p]{1+ m_1^p}\,\right)^{p-1}+
\left (\sqrt[p]{1 + m_1^p}\,\right)^{p-1}=(1+m \, m_1^{p-1})\left(\sqrt[p]{1 + m_2^p}\,\right)^{p-1}.
$$
Consequently:  
$$\begin{array}{r}
  m_2^{p-1} \left(\sqrt[p]{1  + m_1^p }\,\right)^{p-1}=\\ [0.25cm]
\ds\frac{\left(\sqrt[p]{1  + m_2^p }\,\right)^{p-1}+ m\, m_1^{p-1}\left(\sqrt[p]{1  + m_1^p}\,\right)^{p-1}-\left(\sqrt[p]{1  + m_1^p }\,\right)^{p-1}}{m}
\end{array}$$

Therefore, the denominator of the last expression is:
$$\begin{array}{rcl}
m_2^{p-1}\,\left(\sqrt[p]{1  + m_2^p }\,\right)^{p-1}-m_1^{p-1} \left(\sqrt[p]{1  + m_1^p }\,\right)^{p-1}&=& \\[0.45cm]
\frac{\left(\sqrt[p]{1+m_2^p}\,\right)^{p-1}+ m\, m_1^{p-1}\left(\sqrt[p]{1+ m_1^p}\,\right)^{p-1}-\left(\sqrt[p]{1+m_1^p}\,\right)^{p-1}- m\, m_1^{p-1}\left(\sqrt[p]{1 + m_1^p}\,\right)^{p-1}}{m}&=&
\\[0.45cm]
\frac{\left(\sqrt[p]{1  + m_2^p }\,\right)^{p-1}-\left(\sqrt[p]{1  + m_1^p }\,\right)^{p-1}}{m}&&
\end{array}
$$
and hence:
$$\begin{array}{lcl}
-\ds\frac{\frac{\partial{\cal P}_r (a_0,b_0)}{\partial x}}{\frac{\partial {\cal P}_r (a_0,b_0)}{\partial y}}&=&
\frac{\left(\sqrt[p]{{1  + m_2^p }}\,\right)^{p-1}-\left(\sqrt[p]{{1  + m_1^p }}\,\right)^{p-1}}{m_2^{p-1}\,\left(\sqrt[p]{{1  + m_2^p }}\,\right)^{p-1}-m_1^{p-1}\left(\sqrt[p]{{1  + m_1^p }}\,\right)^{p-1}}\\
&=&\frac{\left(\sqrt[p]{{1  + m_2^p }}\,\right)^{p-1}-\left(\sqrt[p]{{1  + m_1^p }}\,\right)^{p-1}}{\frac{\left(\sqrt[p]{{1  + m_2^p }}\,\right)^{p-1}-\left(\sqrt[p]{{1  + m_1^p }}\,\right)^{p-1}}{m}}\\&=&m
\end{array}$$

\item If $P=(x,y) \in \Delta_r$, then the slope $\mu  = \frac{{y - b}}{{x - a}}$  verifies that
$0 < m_2  < \mu  < m_1 .$

{\it Proof: }\\
If $P=(x,y) \in \Delta_r$, then $P$  is placed on the right of $A'$, and hence   $\mu  < m_1 $. 
Moreover, $m_2  < \mu $ since the straight lines $P'P$ and  $AP$ must be concurrent in $P$.\\
Furthermore, the lower bound is $\mathop {\lim }\limits_{x \to \infty } \mu (x) = m_2 $, since, by  rewriting the equation of $\Delta_r$, the following equation is obtained:
$$
A_r(x - a) + B_r(y - b) + A_r \,a + B_r \,b + C_r = (x - a)\sqrt[p]{{1 + \mu ^p }}
$$
and dividing by $(x-a)$, 
$$
A_r + B_ r \,\mu  + \frac{{A_r \,a + B_r \,b + C_r}}{{x - a}} = \sqrt[p]{{1 + \mu ^p }}
$$
Hence 
$$
\mathop {\lim }\limits_{x \to \infty } \left( {A_r + B_r \,\mu  + \frac{{A_r \,a + B_r \,b + C_r}}{{x - a}}} \right) = A_r + B_r \,\mu _\infty   = \sqrt[p]{{1 + \mu _\infty ^p }}$$
Therefore, $$
\frac{{1 + \mu _\infty   \cdot m_2^{p - 1} }}{{\left( {\sqrt[p]{{1 + m_2^p }}} \right)^{p - 1} }} = \sqrt[p]{{1 + \mu _\infty ^p }}$$ with a unique solution (see \cite{Franco0}).
%%%%%%%%%%%%%%%%%%%%%55
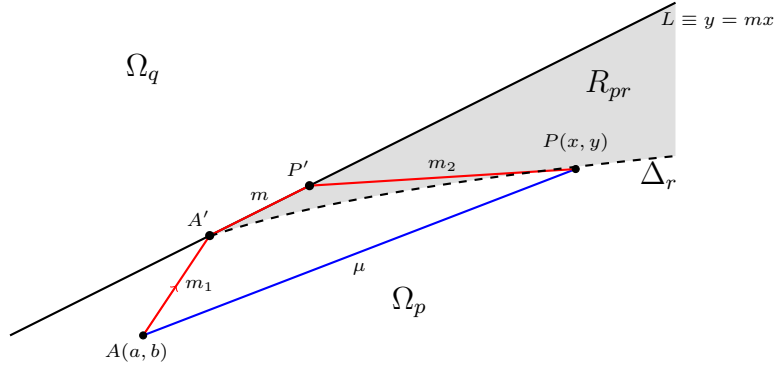
\begin{figure}[ht]\centering
\begin{tikzpicture}[xscale=1.10/1.25,yscale=1.10/1.25]
\fill[color=gray!25]  (10,5)--(3,1.5)-- plot [domain=3:10] (\x, {-0.30+ ln(2*\x) } ) -- cycle;
%\draw[fill=gray!15,ultra thin,dotted] (3,1.5)--(10,1.86)--(10,5)--(3,1.5) ;
\draw[-,thick] (0,0) --(10,5)node[below]{\footnotesize $\phantom{xxxxxxx}L\equiv y=m x$};
%\draw[ultra thin, help lines] (0,0) grid (10,5);
\node at (6,0.5) {\large $\Omega_p$};
\node at (2,4) {\large $\Omega_q$};
\node at (9,3.75) {\large $R_{pr}$};
\node at (1.9,-0.25){\footnotesize $A(a,b)$};
\draw[red,thick] (3,1.5)-- (2,0);
\node at (3.75,2.1){\footnotesize $m$};
\node at (2.85,0.75){\footnotesize $m_1$};
\node at (6.5,2.55){\footnotesize $m_2$};
\node at (5.25,1.0){\footnotesize $\mu$};
\draw[blue,thick] (2,0)--(8.5,2.5);
\draw[red,thick] (3,1.5)-- (4.5,2.25)--(8.5,2.5);
\node at (8.5,2.9){\footnotesize $P(x,y)$};
\draw[->,red] (2,0)--(2.5,.75);
\node at (9.75,2.4){\large $\Delta_r$};
\draw [dashed,thick,domain=3:10] plot(\x, {-0.30+ ln(2*\x) } );    
\filldraw [black] (2,0) circle (1.5pt);
\filldraw [black,thick] (4.5,2.25) circle (1.5pt) node[above] {\footnotesize $P'$\phantom{xx}};
\filldraw [black] (8.5,2.5) circle (1.5pt);
\filldraw [black,thick] (3,1.5) circle (1.5pt) node[above] {\footnotesize $A'$\phantom{xx}};
\end{tikzpicture}
\caption {Graphical representation of $\Delta_r$ and the slopes $m_1$, $m_2$ and $\mu$,}
 \label{fig:slopes}
\end{figure}

\item The implicit function that determines $\Delta_r$ is increasing and concave in  its domain.\\
{\it Proof: } \\
The implicit function of $\Delta_r$ is: $A_r x + B_r y + C_r - \sqrt[p]{{(x - a)^p  + (y - b)^p }}=0$, and
by deriving, it is obtained that:
\begin{equation}\label{ec2}
A_r + B_r \,y' = \frac{{(x - a)^{p - 1}  + (y - b)^{p - 1} y'}}{{\left( {\sqrt[p]{{(x - a)^p  + (y - b)^p }}} \right)^{p - 1} }} = \frac{{1 + \mu ^{p - 1} y'}}{{\left( {\sqrt[p]{{1 + \mu ^p }}} \right)^{p - 1} }},
\end{equation}
and by grouping the terms of the equation, it is obtained that:
$$ A_r - \frac{1}{{\left( {\sqrt[p]{{1 + \mu ^p }}} \right)^{p - 1} }} =
\left(\frac{{\mu ^{p - 1}}}{{\left( {\sqrt[p]{{1 + \mu ^p }}} \right)^{p - 1} }} - B_r\right ) \,y' $$
that is to say:
\begin{equation}\label{Gr}
 G_r = H_r \,y'.
\end{equation}
where  
$$G_r = \frac{1}{{\left( {\sqrt[p]{{1 + m_2^p }}} \right)^{p - 1} }}\; - \frac{1}{{\left( {\sqrt[p]{{1 + \mu ^p }}} \right)^{p - 1} }} > 0$$
and 
$$H_r = \left( {\frac{\mu }{{\sqrt[p]{{1 + \mu ^p }}}}} \right)^{p - 1}  - \left( {\frac{{m_2 }}{{\sqrt[p]{{1 + m_2^p }}}}} \right)^{p - 1}  > 0$$
 since $m_2  < \mu $ and the function $f(x) = \frac{x}{{\sqrt[p]{{1 + x^p }}}}$ is increasing for $x \, > \, 0$. 
 Consequently, $y^{'} \, >0$ and therefore the $p$-curve $\Delta_r$ is increasing.\\

  With respect to the concavity, equation (\ref{ec2}) can be written as
$$A_r + B_r  \,y' = \left( {1 + \mu ^{p - 1} y'} \right)\left( {1 + \mu ^p } \right)^{\frac{1}{p} - 1}$$ 
and by deriving again,  the expression is obtained:
$$\begin{array}{rcl}
B_r \,y'' &=& \left( {(p - 1)\mu ^{p - 2} \mu 'y' + \mu ^{p - 1} y''} \right)\left( {1 + \mu ^p } \right)^{\frac{1}{p} - 1} \\
& &+  \left( {1 + \mu ^{p - 1} y'} \right)\left( {\left( {{\textstyle{1 \over p}} - 1} \right)\left( {1 + \mu ^p } \right)^{\frac{1}{p} - 2} p\mu ^{p - 1} \mu '} \right)
\end{array}
$$
and by simplifying, the following expression is obtained:
$$\begin{array}{rcl}
-H_r \cdot y'' &=& \mu ^{p - 2} \mu '\left( {1 + \mu ^p } \right)^{\frac{1}{p} - 2} (p - 1)(y' - \mu )   \\
& =& \mu ^{p - 2} \left( {1 + \mu ^p } \right)^{\frac{1}{p} - 2} (p - 1)\frac{{(y' - \mu )^2 }}{{x - a}} > 0
  \end{array}
$$
and since $H_r \, > \, 0$, then $y^{''} \, < \, 0$ and hence the curve is concave.
\end{itemize}

%%%%%%%%%%%%%%%%%%%%%%%%%%%%%55
\subsection{Left branch on $\Delta$ ($\Delta_l$)}
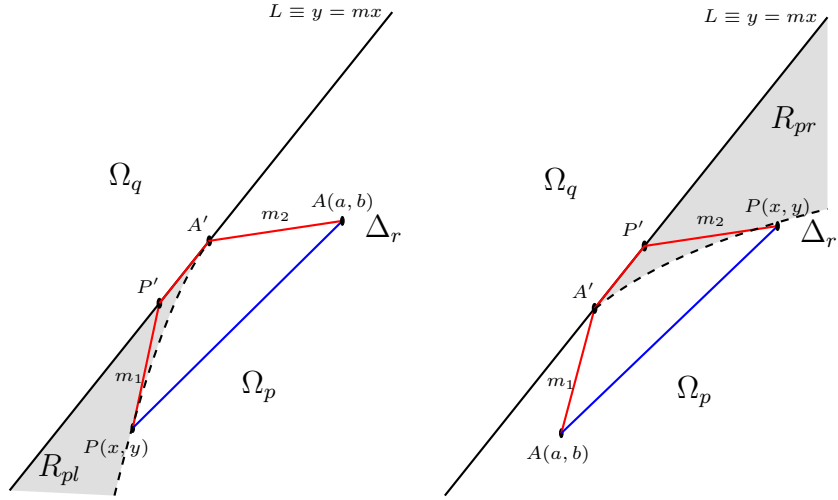
\begin{figure}[ht!]\centering
\begin{tabular}{cc}  
\begin{tikzpicture}[xscale=1.10/2.5,yscale=1.10/1]
\fill[color=gray!25]  plot [domain=-4.10:-1.25] (\x+5.75, -0.2*\x*\x+2.55) --(4.5,2.25)--(-1.5,-0.75)--cycle;   
\draw[-,thick] (-1.5,-0.75)--(0,0) --(10,5)node[left]{\footnotesize $\phantom{xxxxxxx}L\equiv y=m x$};
%\draw[ultra thin, help lines] (0,0) grid (10,5);
\node at (6,0.5) {\large $\Omega_p$};
\node at (2,3) {\large $\Omega_q$};
\node at (0,-.5) {\large $R_{pl}$};
\filldraw [black] (2.2,0) circle (1.5pt);
\node at (1.7,-0.25){\footnotesize $P(x,y)$};
\draw[red,thick] (3,1.5)-- (2.2,0);
\filldraw [black,thick] (3,1.5) circle (1.5pt) node[above] {\footnotesize $P'$\phantom{xx}};
\filldraw [black,thick] (4.5,2.25) circle (1.5pt) node[above] {\footnotesize $A'$\phantom{xx}};
\draw[blue,thick] (2.2,0)--(8.5,2.49);
\draw[red,thick] (3,1.5)-- (4.5,2.25)--(8.5,2.49);
\filldraw [black] (8.5,2.49) circle (1.5pt);
\node[above] at (8.5,2.49){\footnotesize $A(a,b)$};
\node at (2.10,0.60){\footnotesize $m_1$};
\node at (6.5,2.55){\footnotesize $m_2$};
%\draw[->,red] (2.2,0)--(2.5,.75);
\%draw[thick,dashed](4.5,2.25)--(2.8,-0.86);
\draw [dashed, thick,domain=-4.10:-1.25] plot(\x+5.75, -0.2*\x*\x+2.55);
\node at (9.75,2.4){\large $\Delta_r$};
\end{tikzpicture}
&  
\begin{tikzpicture}[xscale=1.10/2.5,yscale=1.10/1]
\fill[color=gray!25]  (10,5)--(3,1.5)-- plot [domain=3:10] (\x, {-0.30+ ln(2*\x) } ) -- cycle;   
\draw[-,thick] (-1.5,-0.75)--(0,0) --(10,5)node[left]{\footnotesize $\phantom{xxxxxxx}L\equiv y=m x$};
%\draw[ultra thin, help lines] (0,0) grid (10,5);
\node at (6,0.5) {\large $\Omega_p$};
\node at (2,3) {\large $\Omega_q$};
\node at (9,3.75) {\large $R_{pr}$};
\filldraw [black] (2,0) circle (1.5pt);
\node at (1.9,-0.25){\footnotesize $A(a,b)$};
\draw[red,thick] (3,1.5)-- (2,0);
\filldraw [black,thick] (3,1.5) circle (1.5pt) node[above] {\footnotesize $A'$\phantom{xx}};
\filldraw [black,thick] (4.5,2.25) circle (1.5pt) node[above] {\footnotesize $P'$\phantom{xx}};
\draw[blue,thick] (2,0)--(8.5,2.49);
\draw[red,thick] (3,1.5)-- (4.5,2.25)--(8.5,2.49);
\filldraw [black] (8.5,2.49) circle (1.5pt);
\node[above] at (8.5,2.49){\footnotesize $P(x,y)$};
\node at (2.00,0.60){\footnotesize $m_1$};
\node at (6.5,2.55){\footnotesize $m_2$};
%\draw[->,red] (2,0)--(2.5,.75);
\draw [dashed,thick,domain=3:10] plot(\x, {-0.30+ ln(2*\x) } );
\node at (9.75,2.4){\large $\Delta_r$};
\end{tikzpicture}
\end{tabular}
\caption {Situation to obtain the Left branch on $\Delta$ (left) and the Right branch on $\Delta$ (right).}
 \label{fig:doslados}
\end{figure}

The developments to obtain the Left branch on $\Delta$ are similar to those seen previously to obtain $\Delta_r$: the key point is to interchange points $P$ and $A$ and slopes $m_1$ and $m_2$, as can be observed in \figurename{ \ref{fig:doslados}}.
%%%%%%%%%%%%%%%%%%%%%
% %%%%%%%%%%%%%%%%%%%%%%%%%%%%%%%%%%%%%%%%%%%%%%%%%%%%%%%%%%%%%%

Thus, the expression of the Left branch of $\Delta$, $\Delta_l$ is as follows:
\begin{equation}\label{ec2l}
\Delta_l  \equiv \left\{ P\in \Omega_p \mbox{ such that }A_l\, x + B_l\, y + C_l - \sqrt[p]{{(x - a)^p  + (y - b)^p }}=0\right\}
\end{equation}
 where %it has been taken into account the equation (\ref{fr2}) obtained from (\ref{fr5}) \cite{Franco1}:\\
 \begin{equation}\label{ec1-3f}
 A_l = \ds\frac{-1}{{\left( {\sqrt[p]{{1 + m_1^p }}} \right)^{p - 1} }},
 B_l =\ds -\left( {\frac{{m_1 }}{{\sqrt[p]{{1 + m_1^p }}}}} \right)^{p - 1}, \mbox{ and }
C_l = \ds  \frac{{a + m_2^{p - 1} b}}{{\left( {\sqrt[p]{{1 + m_2^p }}} \right)^{p - 1} }}.
\end{equation}

\begin{figure}[h]\centering
\begin{tabular}{c}
  \includegraphics[width=60mm]{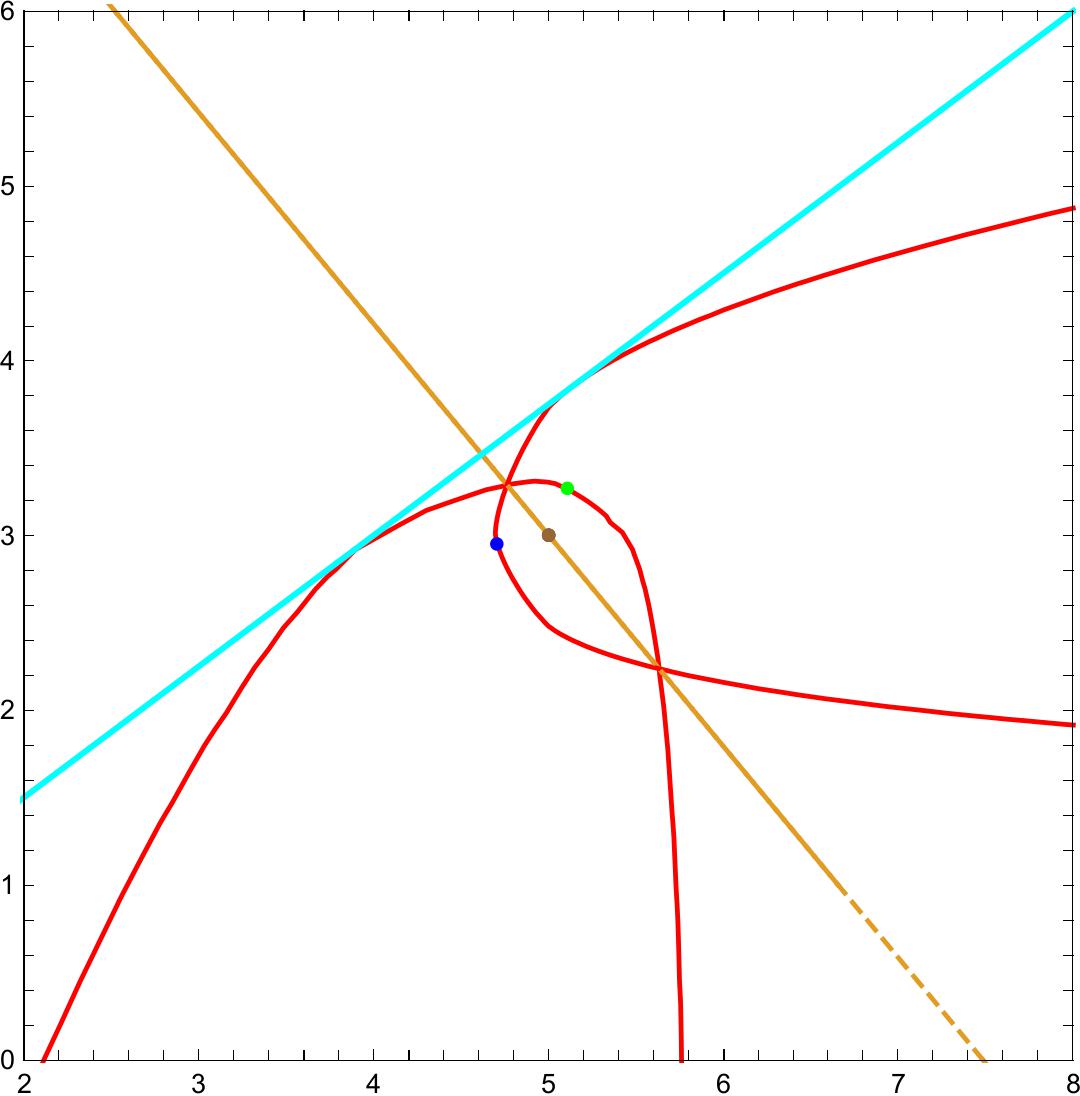} 
 \end{tabular}
\caption {Graphical representation of the $p$-parabolas ${\cal P}_r$ and ${\cal P}_l$ (both in red). The straight $L\equiv y=m\,x$  is depicted in magenta, and the straight line  perpendicular to $L$ is the axis of the both parabolas. The point in blue is the vertex of ${\cal P}_r$ and the point in green is the vertex of ${\cal P}_l$. The point in brown is the focus of both parabolas.} \label{fig:parabolas}
 \end{figure}

 Furthermore, if the $p$-parabola ${\cal P}_l$, of the directrix $L_r\equiv A_l x + B_l y + C_l=0$ and the focus point $A(a,b)$ is consided (see \figurename{ \ref{fig:parabolas}}), then $\Delta_l$ represents, from point $A'=Proj(A, L, m_2)=\left( {\frac{{b - m_2  \cdot a}}{{m - m_2 }}  ,m\frac{{b - m_2  \cdot a}}{{m - m_2 }}} \right)$, the Left branch of the $p$-parabola ${\cal P}_l$.
%%%%%%%%%%%%%%%%%%%%%
%%%%%%%%%%%%%%%%%%%%%%%555

\begin{nota}
From equations (\ref{ec1-3}) and (\ref{ec1-3f}), it holds that
$$ C_r  = A_l \, a+ B_l \, b\mbox{ \phantom{xxxx} and \phantom{xxxx}}
C_l=A_r \, a+ B_r \, b
$$
Therefore:\\
\begin{equation}\label{ec1-6}
\begin{array}{rclcl}
{\cal P}_r & \equiv & A_r \,x + B_r \,y + A_l \, a+ B_l \, b - \sqrt[p]{{(x - a)^p  + ( y-b)^p }}&=&0, \mbox{ and} \\
{\cal P}_l & \equiv & A_l \,x + B_ l \,y + A_r \, a+ B_r \, b - \sqrt[p]{{(x-a)^p  + (y-b)^p }}&=&0
\end{array}
\end{equation}
\end{nota}

The properties of $\Delta_l$ and ${\cal P}_l$ are identical to the properties obtained  by $\Delta_r$ and ${\cal P}_r$ by changing the suffix $r$ with $l$. Thus
\begin{itemize}
 \item $ -1 < A_r , B_r < 0 $
 \item By considering $\bar{p}$, it follows that: $A_l^{\bar{p}}+B_l^{\bar{p}}=1$
\item It is worth noting that if $1<p<q$ and $A \in  \,\Omega _p $ are fixed, then the $\Delta_l $ splits  the half-plane $\Omega _p $, in such a way that:
\begin{itemize}
\item if $P$ is below the half-plane then, $\Delta_l<0 \equiv A_l x + B_l y + C_l < \sqrt[p]{{(x - a)^p  + (y - b)^p }},
$ and the shortest path between point  $A$ and point $P(x,y)$ is a polygonal with three segments with distance $A_l x + B_l y + C_l$; 
\item otherwise, $\Delta_l\geq 0$, and the shortest path between point $A$ and point $P(x,y)$ is the segment joining these points with distance $d_p(A,P)$.  

Thus, let us denote this region as $R_{pl}$ (see the graphical representation on the  left in \figurename{ \ref{fig:doslados}}).
\end{itemize}

 \item The vertex of the the $p$-parabola ${\cal P}_l$ is $\ds
V_r=\left(\frac{a+x_0}{2}\,,\, \frac{b+y_0}{2}\right )$
where
$$
x_0= \ds \frac{a\, B_l^2- A_l (b \, B_l + C_l)}{A_l^2+ B_l^2}, \mbox{ and }
y_0= \frac{b\, A_l^2- B_l (a \, A_l + C_l)}{A_l^2+ B_l^2}
$$

\item Point $A'= Proj(A,L, m_2)=\left( {\frac{{b - m_2  \cdot a}}{{m - m_2 }}  ,m\frac{{b - m_2  \cdot a}}{{m - m_2 }}} \right) $ belongs to  ${\cal P}_l$.

\item The tangent line  to  the $p$-parabola ${\cal P}_l$ at point $A'$ has the slope of the straight line $L\equiv y= m \, x$.

\item The implicit function that determines $\Delta_l$ is increasing and concave in  its domain.
\end{itemize}

%Let us see some consequences of the previous results:

%%%%%%%%%%%%%%%%%%%%%%%%%%%%%%%%%%%%%%%%%%
\subsection{Consequences}
\begin{itemize}
% \item Given the p-parabola $\Delta_r$ and the point $A'= (a_0,b_0)$, we know that $A'\in \Delta_r$ and the %tangent line to $\Delta_r$ at point $A'$ has slope m. If $(x,y)\in \Delta_r$ , for $x>a$, then ${\mu= \frac{y-b}{x-a}}%$ satisfies that $0< m_2<\mu < m_1$ fulfilling that $ \mathop {\lim }\limits_{x \to \infty }\mu=m_2$. $\Delta_{r}$ is %an increasing and concave function for $x>a_0, \, y> b_0$. Therefore the p-parabola remains below the line $L$ %in  $R_{pr}$. In a similar way, the same study is done for the $\Delta_l$ parabola and the point $A''= (a_1,b_1)$ %obtaining the $R_{pl}$ region.

\item Given a point $A=(a,b)$,  two regions $R_{pr}$ and $R_{pl}$ are obtained  from $\Delta_r$ and $\Delta_l$. Let be $R_{pc}=\Omega_p \setminus (R_{pr} \cup R_{pl})$, hence  (see \figurename{ \ref{Rpc01}}):   
$$d(P,A)=\left\{
\begin{array}{rcl}
 A_r\, x+ B_r \, y+C_r& if & P=(x,y)\in R_{pr}\\
 A_l\, x+ B_l \, y+C_l& if & P=(x,y)\in R_{pl}\\
 d_p(A,P)&  if & P=(x,y)\in R_{pc}\\
\end{array}
\right .$$
%%%%%%%%%%%%%%%%%%%%%%55
\begin{figure}[ht]
\centering
\includegraphics[width=80mm]{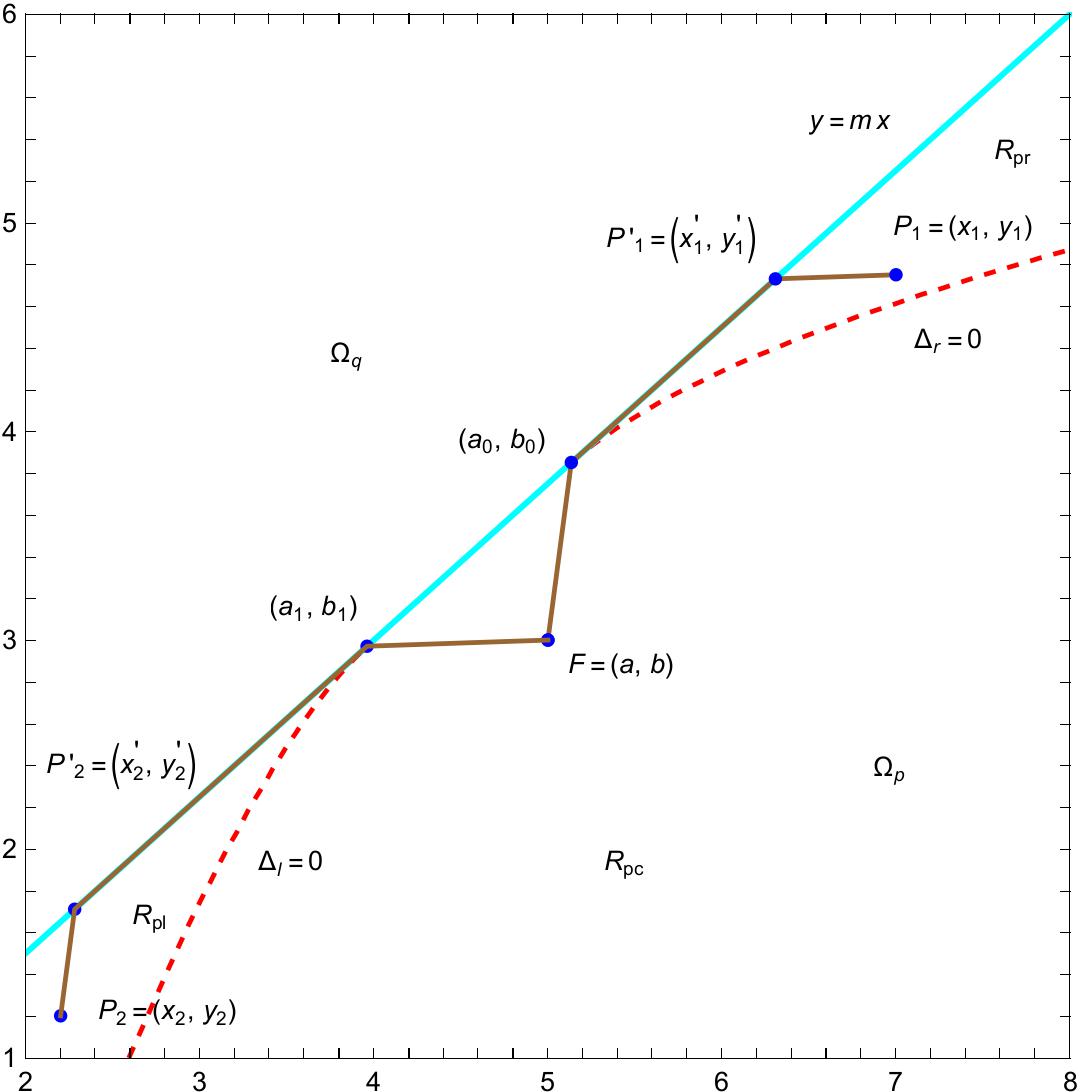}
\caption {$R_{pr}$, $R_{pc}$, $R_{pl}$  Regions in $\Omega_p$.
$Proj(F,L, m_1)=(a_0,b_0)$ and $Proj(F,L, m_2)=(a_1,b_1)$. }
\label{Rpc01}
\end{figure}
%%%%%%%%%%%%%%%%%%%%55   
 
 %\item It is not difficult to proof that $tg(\alpha)= \frac{m+ m_2}{1+ m m_2}$ , $tg(\beta)= \frac{m+ m_1}{1+ m m_1}$ and $tg(\gamma)= \frac{m_1+ m_2}{1+ m_1 m_2}$, where $\alpha $ is the angle formed by the vectors $(v, \, -v_1)$,  $\beta$ the angle formed by $(-v, -v_0)$  and $\gamma$ the angle formed by $v_2= (v_0, v_1)$, with $v_0= (a_0-a, b_0-b)$, $v_1= (a_1-a, b_1-b)$ and $v= (a_0-a_1, b_0-b_1)$.
% 
 \item In the same way as in \cite{Franco2}, two p-curves have been found splitting $\Omega_p$ into three regions delimited by the upper branch of the $\Delta_r=0$ and the Left branch of  $\Delta_l=0$  (see \figurename{ \ref{Rpc01}}).
 \item Note that if $q= \infty$, then $m_2=0$ in equation (\ref{fr5}), $A_l=-1$, $B_r=0$, and $C_l=a$.
\end{itemize}  

%%%%%%%%%%%%%%%%%%%%%%%%%%%5
\section{The Parabola-MFP (P-MFP($p,q$)) Procedure}\label{section_pMFP}
The algorithm to be developed  to solve the problem  (\ref{ec1-2}) is the same that as presented in  \cite{Franco2} in the region $\Omega_q$, since in this region nothing  new is given in this paper. With respect to the region $\Omega_p$, the results on the p-curves are taken into account.

%%%%%%%%%%%%%%%55555
\subsection{Solution of the problem  (\ref{ec1-2}) in $\Omega_p$ }\label{subse:omegap}
Given $X =(x,y)\in \Omega_p$, by using the $D(X)$ function of problem (\ref{ec1-2}) and the previous results (see \figurename{ \ref{Rpc01}}), the following function $D_p(X)$ is considered, in order to solve the location problem (\ref{ec1-2}):
\begin{equation}\label{ec-D}
\begin{array}{rcl}
D_p(x,y)&=& \ds\sum_{i \in S_{pc}} \omega_i\, \sqrt[p]{ ((x-a_i)^2+\epsilon)^{\frac{p}{2}}+((y-b_i)^2+\epsilon)^{\frac{p}{2}}}+\\
&&\ds\sum_{i \in S_{pr}} \omega_i \,(A_r\, a_i+ B_r\, b_i+A_l\, x+B_l\, y)+\\
&& \ds\sum_{i \in S_{pl}} \omega_i\, (A_l\, a_i+ B_l\, b_i+A_r\, x+B_r\, y)+\\
&&\ds\sum_{j \in S_{q}} \omega_j \, \sqrt[p]{((x-z_j)^2+\epsilon)^{\frac{p}{2}}+((y-m\,z_{j})^2+\epsilon)^{\frac{p}{2}}}+\\
&&\ds \sum_{j \in S_{q}} \omega_j \, \sqrt[q]{ ((z_{j}-c_j)^2+\epsilon)^{\frac{q}{2}}+((m \, z_{j}-d_j)^2+\epsilon)^{\frac{q}{2}}}
\end{array}
\end{equation}
where $ S_{pc}= S_p \cap R_{pc}$, $ S_{pl}= S_p \cap R_{pl}$, $ S_{pr}= S_p \setminus (S_{pr}\cup S_{pl})$  and $\{z_j\}_{j \in S_p}$ are the gate points needed in order to find the solution point $X$.

With respect to $D(X)$, $D_p(x,y)$ is a smoothing function using the hyperbolic approximation. Furthermore, $D_p(x,y)$ is convex function (\cite{Brimberg1}) and the weights between pairs of new facilities $z_i's$ are all equal to zero. Hence, taking partial derivatives with respect to $x, \, y$ and $\, z_i's$, similar formulae to the ones for the multi-facility problem (MFP) (\cite{Brimberg0},\cite{Brimberg1}) are obtained. Let us see it: if
$$
\begin{array}{lll}
A_{1 i}^k  =  (x^k- a_i)^2 + \epsilon, & A_{2 i}^k  =  (y^k- b_i)^2 + \epsilon,
& A_i^k  =  (A_{1 i}^k)^{\frac{p}{2}} +(A_{2 i}^k)^{\frac{p}{2}}\\
B_{1 j}^k  =  (x^k- z_j^k)^2 + \epsilon, &B_{2 j}^k  =  (y^k- m \, z_j^k)^2 + \epsilon,
& B_j^k  =  (B_{1 j}^k)^{\frac{p}{2}} +(B_{2 j}^k)^{\frac{p}{2}}\\
C_{1 j}^k  =  (z_j^k-c_j)^2 + \epsilon, & C_{2 j}^k  =  ( m \, z_j^k-d_j)^2 + \epsilon,&
C_j^k  =  (C_{1 j}^k)^{\frac{q}{2}} +(C_{2 j}^k)^{\frac{q}{2}}
\end{array}
$$  
Then $D_p(x,y)$ can be written as:
$$
\begin{array}{rcl}
D_p(x,y)&=& \ds \sum_{i \in S_{pc}} \omega_i\, (A_i^k)^{\frac{1}{p}}+
 \ds\sum_{i \in S_{pr}} \omega_i \,(A_r\, a_i+ B_r\, b_i+A_l\, x+B_l\, y)+\\
&&\ds\sum_{i \in S_{pl}} \omega_i\, (A_l\, a_i+ B_l\, b_i+A_r\, x+B_r\, y)+
\ds\sum_{j \in S_{q}} \omega_j \, (B_j^k)^{\frac{1}{p}}+\\
&&\ds \sum_{i \in S_{q}} \omega_j \, (C_j^k)^{\frac{1}{q}}
\end{array}
$$   
Therefore and in the same way as done by \cite{Brimberg1}:
\begin{equation}\label{xk1}\tiny
x^{k+1}=\ds \frac{\ds\sum_{i \in \, S_{pc}}\!\!\!  w_i  a_i  (A_{1 i}^k)^{\frac{p}{2}-1}  (A_i^k)^{\frac{1}{p}-1} +\!\sum_{j \in \, S_{q}}\!\!\! w_j  z_j^k   (B_{1  j}^k)^{\frac{p}{2}-1} (B_j^k)^{\frac{1}{p}-1} - A_l\! \sum_{i \in S_{pr}} \omega_i- A_r \! \sum_{i \in S_{pl}} \omega_i}{\ds\sum_{i \in \, S_{pc}} w_i  (A_{1  i}^k)^{\frac{p}{2}-1}  (A_i^k)^{\frac{1}{p}-1}+\sum_{j \in \, S_{q}} w_j    (B_{1  j}^k)^{\frac{p}{2}-1}  (B_j^k)^{\frac{1}{p}-1}}
\end{equation}    
\begin{equation}\label{yk1}  \tiny 
y^{k+1}=\ds\frac{\ds\sum_{i \in \, S_{pc}}\!\!\! w_i  b_i (A_{2 \, i}^k)^{\frac{p}{2}-1} (A_i^k)^{\frac{1}{p}-1}\!\! +\!\!\sum_{j \in \, S_{q}}\!\!\! w_j m z_j^k   (B_{2 j}^k)^{\frac{p}{2}-1} \, (B_j^k)^{\frac{1}{p}-1} - B_l\!\! \sum_{i \in S_{pr}} \omega_i- B_r \!\!\sum_{i \in S_{pl}} \omega_i}
{\ds\sum_{i \in \, S_{pc}} w_i   (A_{2  i}^k)^{\frac{p}{2}-1}  (A_i^k)^{\frac{1}{p}-1}+\sum_{j \in S_{q}} w_j    (B_{2 j}^k)^{\frac{p}{2}-1} \, (B_j^k)^{\frac{1}{p}-1}}
\end{equation}
\begin{equation}\label{zk1}\tiny
z_j^{k+1}= \frac{(B_j^k)^{\frac{1}{p}-1} \left[ (B_{1 \, j}^k)^{\frac{p}{2}-1}\,  x^k+(B_{2 \, j}^k)^{\frac{p}{2}-1}\,  m \, y^k\right]+
(C_j^k)^{\frac{1}{q}-1} \left[ (C_{1 \, j}^k)^{\frac{q}{2}-1} \, c_j+(C_{2 \, j}^k)^{\frac{q}{2}-1}\,m\,   d_j\right]}{(B_j^k)^{\frac{1}{p}-1} \left[ (B_{1 \, j}^k)^{\frac{p}{2}-1}+(B_{2 \, j}^k)^{\frac{p}{2}-1}\,  m^2\right]+(C_j^k)^{\frac{1}{q}-1} \left[ (C_{1 \, j}^k)^{\frac{q}{2}-1} \, +(C_{2 \, j}^k)^{\frac{q}{2}-1}\,m^2\right]}
\end{equation}

Thus, these formulae are used to solve problem (\ref{ec1-2}) in an iterative fashion similar to the Weiszfeld procedure. Nevertheless, it is important to indicate that, in the Hessian matrix of the function $D_p(X)$, only the points of the sets $S_q$ and $S_{pc}$ are taken into account, and this matrix is positive definite (\cite{Brimberg0-1}). Note also that if  the $S_{pr}$ and $S_{pl}$  sets are eventually empty, then the (MFP) Brimberg algorithm is obtained (\cite{Brimberg1}).

Let us indicate the iterative procedure in order to attain the solution of the problem (\ref{ec1-2}) in $\Omega_p$:    
  
\begin{description}      
\item[Step 1.]
A fixed positive tolerance, denoted by $\epsilon_{TOL}$, and a number of iterations, denoted by {\it ITER}, are considered.

\item[Step 2.]
Given $p$, $q$, and $m$, the slopes $m_1$ and $m_2$ are obtained. 

\item[Step 3.]
 Let $k=0$ and $D_0= D_p(X^0)$ such that $X_0$ is the starting point. This point $X_0$  is  taken as the solution point, in which case it either belongs to $\Omega_p$, using MFP algorithm (\cite{Brimberg1}), or a point close to solution point  belonging to $\Omega_p$. Another starting point that can be taken into account is the centre of gravity of all the points in $S_p$ with their corresponding weights.
 
\item[Step 4.] 
$k := k+1$ and  a new $X^{k+1}$ is obtained as follows:                                                       
\begin{description}  
\item [a)] For the point  $X^k$, the two p-curves $\Delta_r$, $\Delta_l$, the points  $(a_0,b_0)=Proj(X^k, L, m_2)$ and  $(a_1,b_1)=Proj(X^k, L, m_1)$, and  the sets $R_{pr}$, $R_{pc}$  and $R_{pl}$ are obtained (Section \ref{sec:pcurve}).  \\
  In this iteration, the points of $\{P_i\}_i\in S_p$ are assigned to each of the corresponding sets $S_{pr}$, $S_{pc}$  and $S_{pl}$.  Subsequently, $x^{k+1}$, $ y^{k+1}$  and the gates $\{z_j^{k+1}\}_{j\in S_p} $ (equations (\ref{xk1}), (\ref{yk1}), and (\ref{zk1})) are obtained and, hence, the new point  $X^{k+1}$   and the corresponding $D_{k+1}=D_p(X^{k+1})$ are obtained.
  
\item [c)] If $|D_{k+1}-D_k |< \epsilon_{TOL}$ or {\it ITER} $\geq k$, then go to d). Otherwise go to a).
\item[d)] $X^*_p= X^{k+1}$ is the solution obtained.
\end{description}
\end{description}   
 
% It is worth noting that $\{z_j^{k+1}\}_{j\in S_p} $ can be obtained solving the problem (\ref{ec1-1}), where $P=X^k$, $X= z_j^{k+1}$, $Q=Q_j$ and then another algorithm, similar to one in \cite{Franco2}, can be used.

Let us study an example of the results of this algorithm. 

%%%%%%%%%%%%%%%%%%%%%%%%%55
\begin{example}\label{ex:01}

Let the points $$\{P_i\}_{i=1}^8=\left\{(0.5, 0.4), (3,2), (4,3), (4,6), (5,5), (5,7), (6,4), (9.5,9)\right\}\in \R^2,$$ $m=0.75$, $p=2$,  and $q=3$. Let us suppose that the solution of problem (\ref{ec1-2}) is in $\Omega_p$.\\
  
In this example (see \figurename {  \ref{fig:Algoritmo}}), the starting point $X_0$ is chosen as the solution of the MFP algorithm and the solution point $X^*$ of the problem (\ref{ec1-2}) is obtained. Hence, $\{P_4 ,\,  P_6\}\in S_q$ , $P_1 \in S_{pl}$, $P_8 \in S_{pr}$, and $\{P_2, P_3, P_5, P_7\} \in S_{pc}$.\\    
Points $\{P_4 ,\,  P_6\}\in S_q$ search for the corresponding Gate point $\{z_j^{k+1}\}_{j\in S_p} $ with distance $\ell_q$ on the line ($y = 0.75\, x$)   and then go to the solution point $X^*$ with distance $\ell_p$ (\figurename { \ref{fig:Algoritmo}}).\\  
Point $P_1$ goes to the line ($y = 0.75\, x$) following direction $m_1$ with distance $\ell_p$ to point $P_1'$, continues along the line to point $Proj(X^+, L, m_2)$ with distance $\ell_q$, and then it goes to the solution point $X^*$ with distance $\ell_p$.\\      
The process is similar for point $P_8$, but in this case the point $Proj(X^+, L, m_1)$ is used.\\  
Points in $S_{pc}$ go to solution point $X^*$ with distance $\ell_p$.
%%%%%%%%%%%%55
\begin{figure}[ht]\centering
\includegraphics[width=80mm]{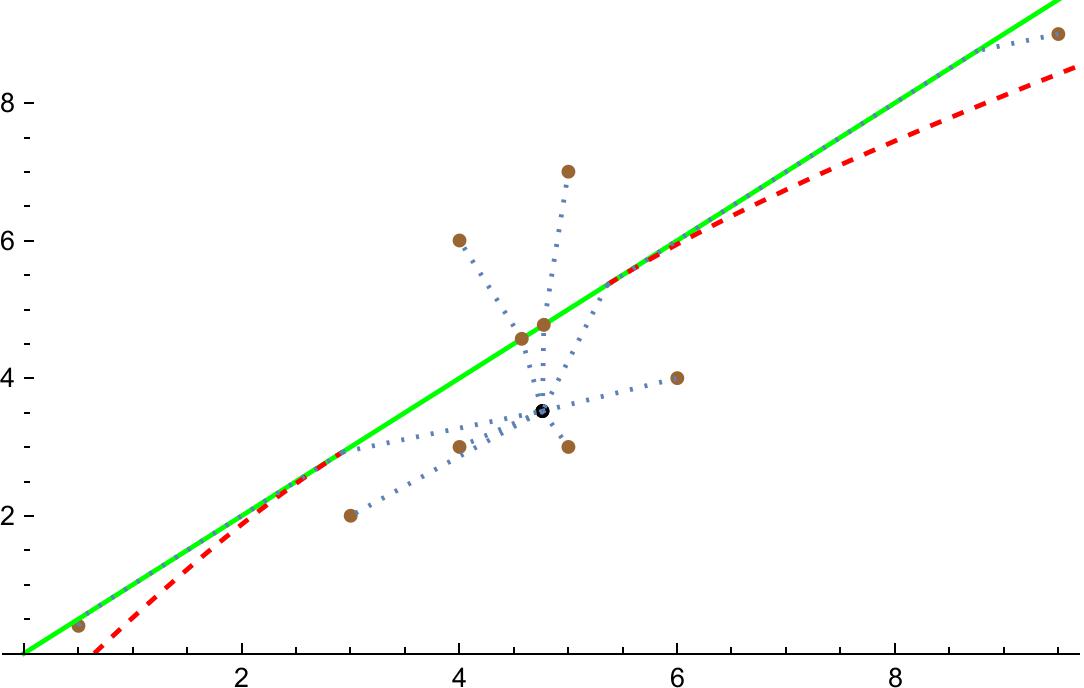}
\caption {Solution of the example \ref{ex:01}. The dashed red  curves  represent the two  p-curves and the point in black is the solution.}
\label{fig:Algoritmo}
\end{figure}
\end{example}
%%%%%%%%%%%%%%%%%%%%55

\begin{figure}[htb!]\centering
\includegraphics[width=95mm]{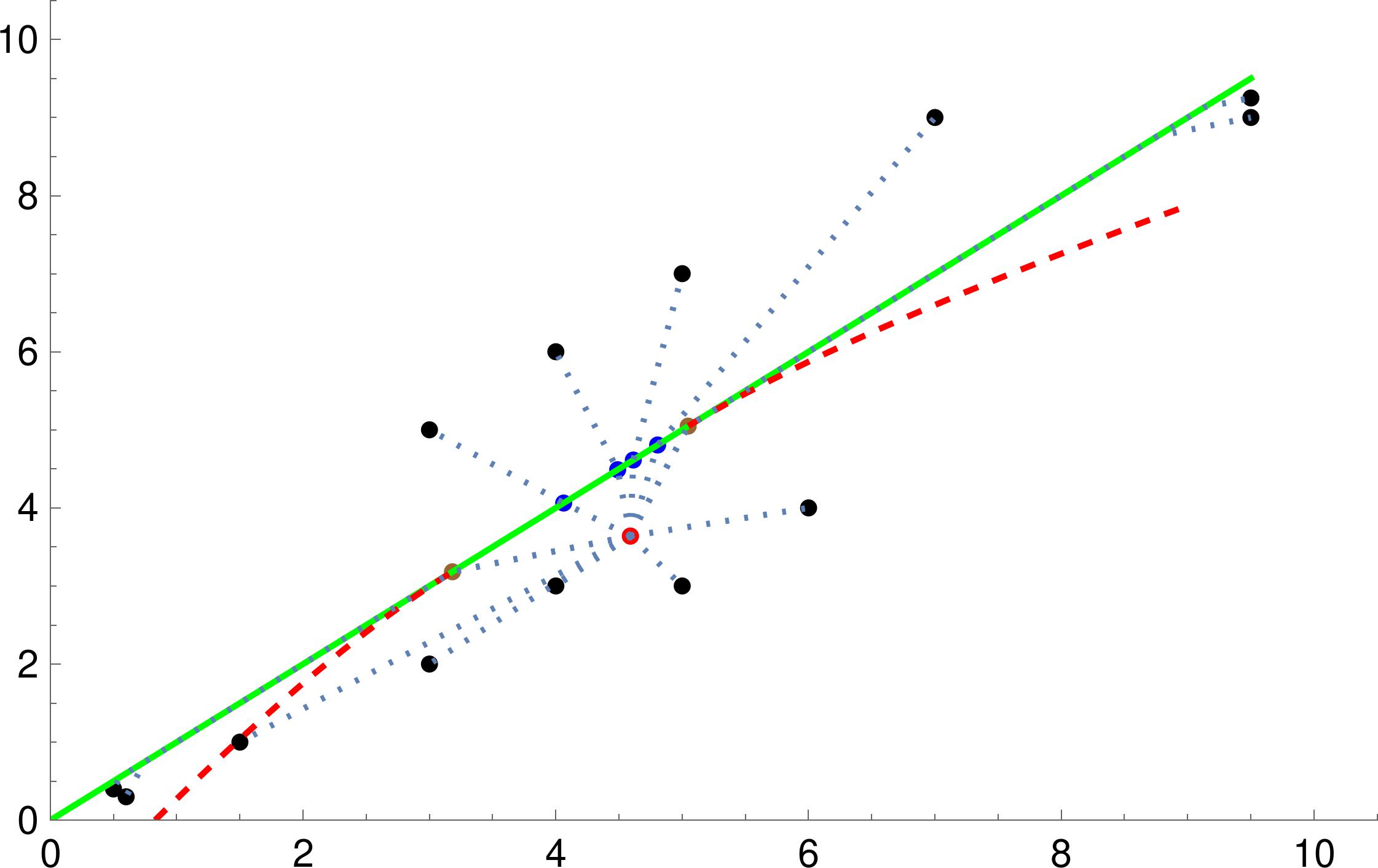}
\caption {Shortest path in Region $\Omega_p$.}
\label{Algoritmo-1}
\end{figure}

\begin{figure}[htb!]\centering
\includegraphics[width=95mm]{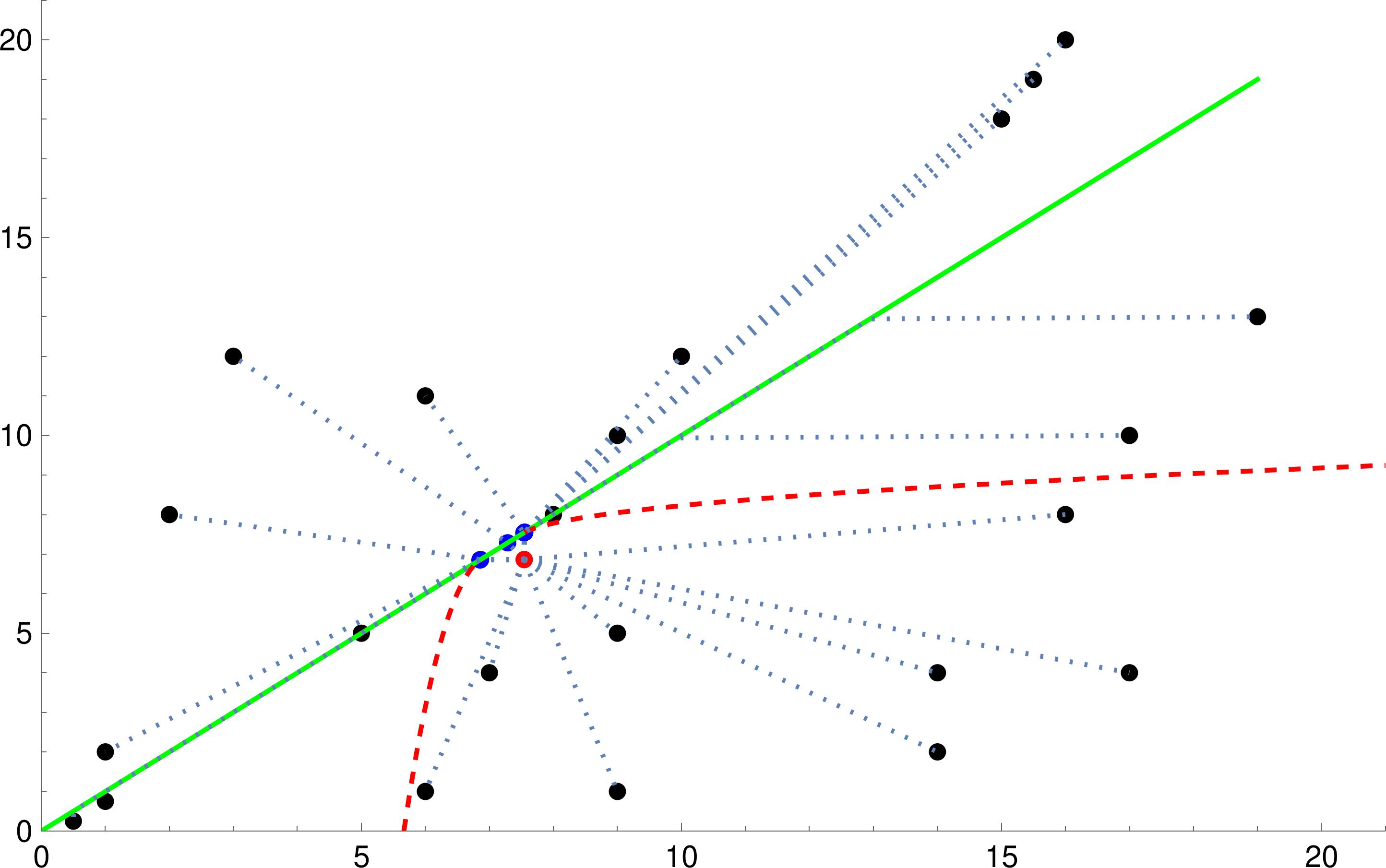}
\caption {Shortest path in Region $\Omega_p$.}
\label{Algoritmo-2}
\end{figure}

\begin{figure}[htb!]
\centering
\includegraphics[width=95mm]{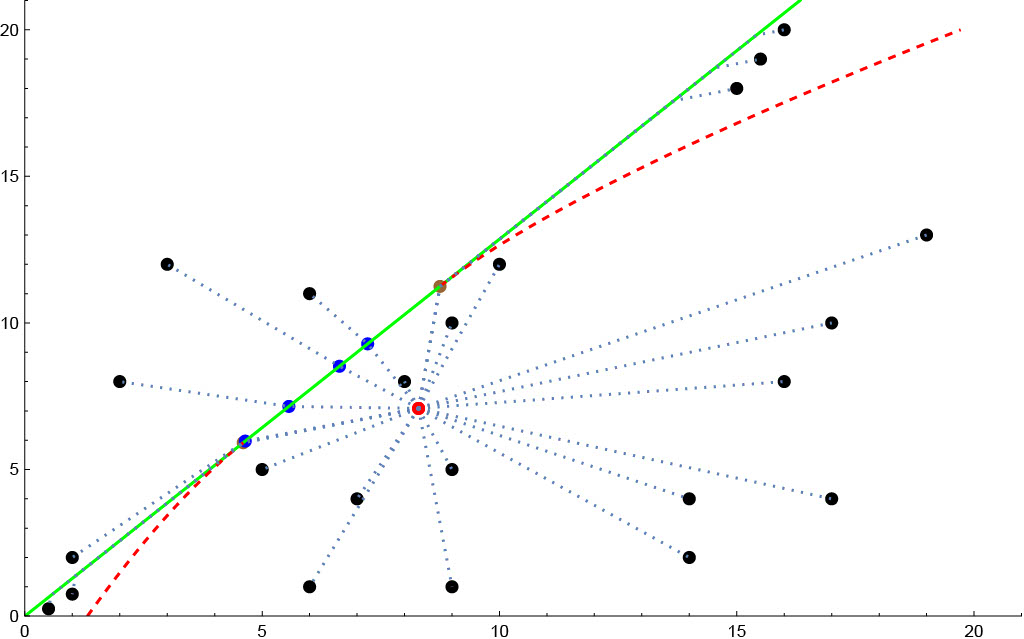}
\caption {Shortest path in Region $\Omega_p$.}
\label{Algoritmo-3}
\end{figure}

Figures \ref{Algoritmo-1}, \ref{Algoritmo-2}, and \ref{Algoritmo-3} show other examples with different data in different regions depending on $p$ and $q$.   
%%%%%%%%%%%%%%%%%5555

%%%%%%%%%%%%%%%%%%%%%%%55555
  
%%%%%%%%%%%%%%%%%%%%%%%%%%%%%%%%%%%
\subsection{Algorithm} \label{Al}
Let us show the P-MFP($p,q$) algorithm:

Let $\cal{P}=$ $\{P_i\}_{i\in S_p}\subset \Omega_p$ and  $\cal{Q}=$ $\{Q_j\}_{j\in S_q}\subset \Omega_q$ be two sets of points with $S_p$ and $S_q$ as two index sets, and  $w_i, w_j \in \R^{+}$ for any $i \in S_p,\, j\in S_q$, where $1<p < q < \infty$ and $L \in \Omega_q$. In order to solve the optimisation problem (\ref{ec1-2}), an algorithm called P-MFP($p$,$q$) is considered as follows:\\
First, the function 
$$D(X)=\,\left\{  \sum_{i\in S_p} w_i \, d(P_i,X)+
\sum_{j\in S_q} w_j \, d_q(Q_j,X)\right\}$$ is considered and the
solution of the three optimisation problems is as follows:\\
\begin{description}    
\item[P1.-] The solution of the problem (\ref{ec1-2}) is sought on the straight line $L$, that is, 
\begin{equation}\label{ec1-5}
\min_{X \in L} \, D(X)
\end{equation}
by using the algorithm given in  \cite{Franco2}. This point solution is denoted by $X^*_L$.

\item[P2.-] The solution of the problem (\ref{ec1-2}) is sought on the region $\Omega_q\setminus L$, that is, 
\begin{equation}\label{ec1-6}
\min_{X \in \Omega_q\setminus L} \, D(X)
\end{equation}
by using the algorithm given in  \cite{Franco2}. This point solution is denoted by $X^*_q$. 
    
\item[P3.-] The solution of the problem (\ref{ec1-2}) is sought in the region $\Omega_p$, that is,
\begin{equation}\label{ec1-7}
\min_{X \in \Omega_p} \, D(X)    
\end{equation}
by using the algorithm given in Section \ref{subse:omegap}. This point solution is denoted by $X^*_p$.    

\item[P4.-] The solution to problem (\ref{ec1-2}) is obtained as follows:
$$X^* = argmin \left \{ D(X^*_L), D(X^*_q), D(X^*_p) \right \}$$
\end{description}

%%%%%%%%%%%%%%%%%%%%%%%%%%55
\section{Experimentation}\label{section6}
The algorithms have been implemented in MATHEMATICA using an Intel(R)
Core(TM)i7-8700 CPU @3.20 GHZ 3.19 GHZ, 8 GB RAM. Tolerance is
taken as $\epsilon_{TOL}=10^{-6}$, the hyperbolic distance $\epsilon =10^{-5}$ and the maximum of iterations is taken as $ITER=10^{4}$.

In order to show the performance of the P-MFP algorithm with respect to the MFP algorithm, two experiments have been carried out. 

%%%%%%%%%%%%%%%%%%%%%%%%%%%%%%%%%%%%%%%555
\subsection{Experiment 1}
Let us consider the example given in \cite{Brimberg1}, where $50$ points are considered. We have considered the values $m=0.5$, $p=1.25$, and $q=2$. Hence, \tablename{ \ref{tab:ejemplo1}} is obtained where a comparative between the solutions given with the two algorithms is provided.
%%%%%%%%%%%%%%%%%%%%%%55
\begin{table}[ht]
%\begin{center}
\caption{Results of the MFP and P-MFP algorithms with the data given in \cite{Brimberg1}.}
 \label{tab:ejemplo1}
 \begin{tabular}{ccc}
 \toprule
  Algorithm &Solution $X^*$ & Distance $D(X^*)$\\
  \midrule
  MFP   & $(10.9993, 8.5146)$ & $1028.34$\\
  P-MFP & $(11.0017, 8.4384)$ & $1027.28$\\
  \botrule  
 \end{tabular}
 %\end{center}
\end{table}
%%%%%%%%%%%%%%%%%%%%%%%%%55

The differences between the distances given in \tablename{ \ref{tab:ejemplo1}} are due to the fact that there are points at which $\Delta_r<0$ or $\Delta_l<0$ where the focus of the $p$-parabolas is the optimal point obtained from the MFP algorithm. 

Therefore, the condition for the P-MFP algorithm to improve the  MFP algorithm is that at least one of the sets $S_{pl}$ and $S_{pc}$ associated to the optimal point obtained from the MFP algorithm is non-empty.

%%%%%%%%%%%%%%%%%%%%%%%%%%%%%%%%%%%%%%%555
\subsection{Experiment 2}
In order to confirm these results, an experiment is carried out by considering   the following values: $p\in\{1.5,1.75\}$, $m\in\{1,1.5\}$, and  $q=2$.

With respect to the points, these are generated uniformly in the square $[0,1]\times[0,1]$. Hence, 100 datasets with different size,  $n\in\{30,50\}$), have been considered in which the optimal MFP solution belongs to $\Omega_p$. 

The result of the experiment can be observed in \tablename{ \ref{tab:simulacion}}. 
As mentioned earlier, in no case is the solution obtained with the P-MFP algorithm worse than that obtained with the MFP: the solution is always equal or better.

Let us discuss the items provided in \tablename{ \ref{tab:simulacion}}. 
\begin{itemize}
 \item The `Improvement(\%)' row gives the percentage of cases in which the P-MFP solution improves the solution obtained with MFP. 
 
 Thus, value 72 should be taken to mean that 72\% of the time the solution is better and in the remaining 28\%, the same solution is  reached.
 
 \item The `Dist. Improv.' row gives the average of the differences between the distances from the MFP algorithm and the P-MFP algorithm.
 
As can be observed, it is positive since the value $D(X)$ will always be greater or equal for the solution given by the MFP algorithm than for the P-MFP algorithm. 

The differences in absolute terms between the solutions are small. However,  the units in which distances are measured and the number of points in $S_{pl}$ and $S_{pr}$'$x$ must be taken into account.

 \item The `Dist. Improv. (\%)' row  provides the average of the differences in percentage terms of the value obtained with the MFP method.
 
 \item The `Points Distance' row gives the average of the $p$-distance between the optimal points obtained through the two procedures.
 
 \item The rows `Points in $S_{pl}$' and `Points in $S_{pr}$' represent the average number of points belonging to those regions. Choosing points at random,  there are sometimes points in $S_{pl}$ and/or $S_{pc}$, but this number of points is small. Therefore,  in all but one case, the average number of points in $S_{pl}$ and/or $S_{pr}$' is less than one. 
\end{itemize}

%%%%%%%%%%%%%%%%%%%%%%%%%%%%%555
\begin{table}[ht]
\caption{Results of the simulation procedure for q=2} 
\label{tab:simulacion}%
%\centering % used for centering table
%\resizebox{\columnwidth}{!}{
\begin{tabular}{lc|c|c|c||c|c|c|c}
\toprule
 & \multicolumn{4}{c||}{$p=1.5$} & \multicolumn{4}{c}{$p=1.75$}\\  %\cline{2-9}
 \midrule
\multicolumn{1}{ c | }{$m$} & \multicolumn{2}{c}{$1.0$}&\multicolumn{2}{c||}{$1.5$}&\multicolumn{2}{c|}{$1.0$}&\multicolumn{2}{c}{$1.5$} \\
\midrule
 \multicolumn{1}{ c | }{$n$} &$30$&$50$&$30$&$50$&$30$&$50$&$30$&$50$ \\
\midrule
\multicolumn{1}{ l | }{Improvement(\%)}&72&91&30&53&52&71&10&16 \\
\midrule
\multicolumn{1}{ l | }{Dist. Improv.}&0.0229&0.0409&0.0056&0.0119&0.0051&0.0099&0.0021&0.0040 \\

\multicolumn{1}{ l | }{Dist. Improv.(\%)}&0.1963&0.2091&0.0473&0.0588&0.0446&0.0519&0.0192&0.0204 \\

\multicolumn{1}{ l | }{Points Distance}&0.0084&0.0088&0.0026&0.0032&0.0035&0.0042&0.0014&0.0014 \\

\multicolumn{1}{ l | }{Points in $S_{pl}$}&0.6900&1.1500&0.3400&0.6700&0.4000&0.6600&0.1300&0.2500\\

\multicolumn{1}{ l | }{Points in $S_{pr}$}&0.6600&0.9600&0.0900&0.2400&0.3900&0.8600&0.1300&0.1900 \\
\botrule
\end{tabular}%}
\end{table}
%%%%%%%%%%%%%%%%%%%%%%555

%%%%%%%%%%%%%%%%%%%%%%%%%%%%%%%%%%%%%%%%%%%%%%%55555555
\section{Conclusions}
The P-MFP algorithm has been presented in order to
solve the optimisation problem (\ref{ec1-2}). This problem was solved in \cite{Franco2} with $p=1$ and $p < q$ and obtained two straight lines by splitting $\Omega_1$ into three regions $\Omega_{1l}$, $\Omega_{1c}$, and $\Omega_{1r}$. In this paper, two $p$-parabolas are found by splitting the $\Omega_p$ region into three regions $\Omega_{pl}$, $\Omega_{pc}$, and $\Omega_{pr}$. 

The MFP algorithm disregards these regions, and therefore fails to take into account that there are two regions ($\Omega_{pl}$ and $\Omega_{pr}$) in which the points therein use a three-link path to reach the solution point. 
This three-link path reduces the distance between the points of these two regions to the solution point. 

Hence, the calculation time is longer in the P-MFP algorithm than it is in the MFP algorithm since, in each iteration, and given the point obtained as a possible solution point, the new regions given by the new $p$-parabolas must be obtained. However, these $p$-parabolas are obtained analytically, as is the location of the points in each of the subregions. Several properties are found describing the properties of $p$-parabolas.

An empirical study has been carried out where it has been proved that the
P-MFP algorithm is more efficient than the MFP algorithm since the solution
provided by the P-MFP algorithm attains a smaller value of the objective function than that given by the MFP algorithm.

%%%%%%%%%%%%%%%%%%%%%%%%%%%%%%%%%%%%%%%%
\section*{Declarations}

\begin{itemize}
%\item Funding: 
\item Conflict of interest: The authors have no conflict of interests to declare that are relevant to the content of this
paper.
\item Ethics approval and consent to participate: `Not applicable'
\item Consent for publication: `Not applicable'
\item Data availability: `Not applicable'
\item Materials availability: `Not applicable'
\item Code availability: `Not applicable' 
\item Author contribution: These authors contributed equally to this work.
\end{itemize}

%%%%%%%%%%%%%%%%%%%%%%%%%%%%%
\bmhead{Acknowledgements}
This research is partially supported by the project of the
Ministerio de Ciencia e Innovaci\'{o}n (PID2021-122208OB-I00) and the project of the
Junta de Andaluc\'{\i}a, Consejer\'{\i}a de Universidad, Investigaci\'{o}n e Innovaci\'{o}n
(ProyExcel-00286).

%%%%%%%%%%%%Referencias%%%%%%%%%%%%%
%\bibliographystyle{elsarticle-harv}
%\bibliography{paco}
%\bibliography{paco}% common bib file

\end{document}